\documentclass{elsarticle}
\usepackage{cite}
\usepackage{graphicx}
\usepackage{color}
\usepackage{verbatim}
\usepackage{url}
\usepackage[numbers]{natbib}

\usepackage{subfigure}

\usepackage{amsmath}
\usepackage{amsfonts}
\usepackage{amssymb}
\usepackage{stmaryrd}

\usepackage{tikz}
\usetikzlibrary{positioning}
\usetikzlibrary{calc}
\usetikzlibrary{chains}

\usepackage[noend]{algorithmic}
\usepackage{algorithm}
\algsetup{indent=2em}

\usepackage[margin=10pt,font=small,labelfont=bf,labelsep=period]{caption}

\newcommand{\assign}{:=}
\newcommand{\BARRIER}{\STATE \rule[0.3em]{2cm}{0.5pt} Memory Fence \rule[0.3em]{2cm}{0.5pt}}
\newcommand{\D}{\mathsf{D}}
\newcommand{\halfopen}[2]{[#1,#2\rangle}
\newcommand{\intceil}[2]{\lceil#1\rceil_{#2}}

\newcommand{\reg}[1]{\mathtt{#1}}
\newcommand{\smem}[1]{\mathtt{#1}^{\mathtt{S}}}
\newcommand{\gmem}[1]{\mathtt{#1}^{\mathtt{G}}}
\newcommand{\tex}[1]{\mathtt{#1}^{\mathtt{T}}}

\newcommand{\mathd}{\,\mathrm{d}}

\newenvironment{myfigure}{\begin{figure}\centering}{\end{figure}}

\newcommand{\avg}[1]{\{#1\}}
\newcommand{\jump}[1]{\left\llbracket#1\right\rrbracket}
\let\cite=\citep

\colorlet{darkgreen}{green!50!black}

\begin{document}

\begin{frontmatter}

\title{Nodal Discontinuous~Galerkin Methods on Graphics~Processors}
\author[brown]{A.~Kl\"ockner}
\ead{andreas\_kloeckner@brown.edu}
\author[rice]{T.~Warburton}
\ead{timwar@rice.edu}
\author[rice]{J.~Bridge}
\ead{jab3@rice.edu}
\author[brown]{J.~S.~Hesthaven}
\ead{jan\_hesthaven@brown.edu}

\address[brown]{Division of Applied Mathematics, Brown University, Providence, RI 02912}
\address[rice]{Department of Computational and Applied Mathematics, Rice University, Houston, TX 77005}

\begin{abstract}

Discontinuous Galerkin (DG) methods for the numerical solution of partial
differential equations have enjoyed considerable success because they are both
flexible and robust: They allow arbitrary unstructured geometries and easy
control of accuracy without compromising simulation stability.  Lately, another
property of DG has been growing in importance: The majority of a DG operator is
applied in an element-local way, with weak penalty-based element-to-element
coupling.

The resulting locality in memory access is one of the factors that enables DG
to run on off-the-shelf, massively parallel graphics processors (GPUs).  In
addition, DG's high-order nature lets it require fewer data points per
represented wavelength and hence fewer memory accesses, in exchange for higher
arithmetic intensity. Both of these factors work significantly in favor of a
GPU implementation of DG.

Using a single US\$400 Nvidia GTX~280 GPU, we accelerate a solver for Maxwell's
equations on a general 3D unstructured grid by a factor of 40 to 60 relative to
a serial computation on a current-generation CPU. In many cases, our algorithms
exhibit full use of the device's available memory bandwidth. Example
computations achieve and surpass 200 gigaflops/s of net application-level
floating point work. 

In this article, we describe and derive the techniques used to reach this
level of performance. In addition, we present comprehensive data on the 
accuracy and runtime behavior of the method.

\end{abstract}


\end{frontmatter}

\section{Introduction}

Discontinuous Galerkin methods
\cite{reed_triangular_1973,cockburn_runge-kutta_1990,hesthaven_nodal_2007} are,
at first glance, a rather curious combination of ideas from Finite-Volume and
Spectral Element methods. Up close, they are very much high-order methods by design.
But instead of perpetuating the order increase like conventional global
methods, at a certain level of detail, they switch over to a decomposition into
computational elements and couple these elements using Finite-Volume-like surface Riemann
solvers. This hybrid, dual-layer design allows DG to combine advantages from
both of its ancestors.  But it adds a third advantage: By adding a movable
boundary between its two halves, it gives implementers an added degree of
flexibility when bringing it onto computing hardware.

A momentous change in the world of computing is now opening an opportunity to
exploit this flexibility even further. Previously, the execution time of a
given code could be determined simply by counting how many floating point
operations it executes.  More recently, memory bottlenecks, in the form of
bandwidth limitation and fetch latency, have taken over as the dominant factors,
and CPU manufacturers use large amounts of silicon to mitigate this effect. It
is quite instructive and somewhat depressing to compare the chip area used for
caches, prediction, and speculation in recent CPUs to the area taken up by the
actual functional units. The picture is changing, however, and graphics
processors, having recently turned into general-purpose programmable units,
were the first to do away with expensive caches and combat latency by massive
parallelism instead.  In this article, we explore how and with what benefit DG
can be brought onto GPUs.

Two main questions arise in this endeavor: First, how shall the computational
work be partitioned? In a distributed-memory setting, the answer is quite
naturally domain decomposition. For the shared-memory parallelism of a GPU,
there are several possibilities, and there is often no single answer that works
well in all settings. Second, DG implementations on serial processors often
rely heavily on the availability of off-the-shelf, pre-tuned linear algebra and
communication primitives. These aids are either unavailable or unsuitable on a
GPU platform, and in stark contrast to the relatively straightforward
implementation of DG on serial machines, optimal use of graphics hardware for
DG presents the implementer with a staggering number of choices. We will
describe these choices as well as a generative approach that exploits them to
adapt the method to both the problem and the hardware at run time.

Using graphics processors for computational tasks is by no means a new idea.
In fact, even in the days of marginally programmable fixed-function hardware,
some (especially particle-based) methods obtained large speedups from running
on early GPUs.  (e.g. \cite{li_implementing_2003}) In the domain of solvers
for partial differential equations, Finite-Difference Time-Domain (FDTD)
methods are a natural fit to graphics processors and obtained speedups of
about an order of magnitude with relative ease (e.g.,
\cite{krakiwsky_acceleration_2004}). Finite Element solvers were also brought
onto GPUs relatively early on (e.g., \cite{goeddeke_accelerating_2005}), but
often failed to reach the same impressive speed gains observed for the simpler
FD methods. In the last few years, high-level abstractions such as Brook and
Brook for GPUs \cite{buck_brook_2004} have enabled more and more complex
computations on streaming hardware. Building on this work, Barth et al.\
\cite{barth_streaming_2005} already predicted promising performance for
two-dimensional DG on a simulation of a the Stanford Merrimac streaming
architecture \cite{dally_merrimac_2003}.  Nowadays, compute abstractions are
becoming less encumbered by their graphics heritage
\cite{lindholm_nvidia_2008, cuda_prog_2008}. This has helped bring algorithms
of even higher complexity onto the GPU (e.g.  \cite{gumerov_fast_2008}).
Taking advantage of these recent advances, this paper presents, to the best of
our knowledge, one of the first general finite-element based solvers that
achieves more than an order of magnitude of speedup on a single real-world
consumer graphics processor when compared to a CPU implementation of the same
method.

A sizable part of this speedup is owed to our use of high-order approximations.
High-order methods require more work per degree of freedom than low-order
methods. This increased arithmetic intensity shifts the method from being
limited by memory bandwidth towards being limited by compute bandwidth.  The
relative abundance of cheap computing power on a GPU makes high-order methods
especially beneficial there.

In this article, we will discuss the numerical solution of linear hyperbolic
systems of conservation laws using DG methods on the GPU. Important examples of
this class of partial differential equations (\emph{PDE}s) include the
second-order wave equation, Maxwell's equations, and many relationships in
acoustics and linear elasticity. Certain nontrivial adjustments to the
discontinuous Galerkin method become necessary when treating nonlinear problems
(see, e.g., \cite[Chapter 5]{hesthaven_nodal_2007}). We leave a detailed
investigation of the solution of nonlinear systems of conservation laws using
DG on a GPU for a future publication, where we will also examine the benefit of
GPU-DG for different classes of PDEs, such as elliptic and parabolic problems.

We will further focus on tetrahedra as the basic discretization element for a
number of reasons. First, it is undisputed that three-dimensional calculations
are in many cases both more practically relevant and more plagued by
performance worries than their lower-dimensional counterparts. Second, they
have the most mature meshing machinery available of all commonly used element
shapes.  And third, when compared with tensor product elements, tetrahedral DG
is both more arithmetically intense and requires fewer memory fetches. Overall,
it is conceivable that tetrahedral DG will benefit more from being carried out
on a GPU.

This article describes the mapping of DG methods onto the Nvidia CUDA
programming model. Hardware implementations of CUDA are available in the form
of consumer graphics cards as well as specialized compute hardware. In
addition, the CUDA model has been mapped onto multicore CPUs with good success
\cite{stratton_mcuda_2008}. Rather than claim an artificial `generality', we
will describe our approach firmly in the context of this model of computation.
While that makes this work vendor-specific, we believe that most of the ideas
presented herein can be reused either identically or with mild modifications to
adapt the method to other, related architectures.  To reinforce this point, we
remark that the the emerging OpenCL industry standard \cite{opencl_2008}
specifies a model of parallel computation that is a very close relative of
CUDA.

The paper is organized as follows: We give a brief overview of the theory and
serial implementation of DG in Section \ref{sec:dg-overview}. The CUDA
programming model is described in Section \ref{sec:gpu-hardware}. Section
\ref{sec:dg-gpu-design} explains the basic design choices behind our approach,
while Section \ref{sec:dg-gpu-implementation} gives detailed implementation
advice and pseudocode. Section \ref{sec:experiments} characterizes our
computational results in terms of speed and accuracy. Finally, in Section
\ref{sec:conclusions} we conclude with a few remarks and ideas for future work.

\section{Overview of the Discontinuous Galerkin Method}

\label{sec:dg-overview}

We are looking to approximate the solution of a hyperbolic system of
conservation laws
\begin{equation} 
  u_t + \nabla \cdot F (u) = 0 \label{eq:claw} 
\end{equation}
on a domain $\Omega = \biguplus_{k=1}^K \D_k \subset \mathbb R^d$ consisting of
disjoint, face-conforming tetrahedra $\D_k$ with boundary conditions
\[ 
  u|_{\Gamma_i} (x, t) = g_i (u (x, t), x, t), \hspace{2em} i = 1, \ldots, b,
\]
at inflow boundaries $\biguplus \Gamma_i \subseteq \partial \Omega$. 
As stated, we will assume the flux function $F$ to be linear.  
We find a weak form of (\ref{eq:claw}) on each element $\D_k$:
\begin{align*} 
  0 & = \int_{\D_k} u_t \varphi + [\nabla \cdot F (u)] \varphi \mathd x\\ 
  & = \int_{\D_k} u_t \varphi - F (u) \cdot \nabla \varphi \mathd x
  + \int_{\partial \D_k} ( \hat{n} \cdot F)^{\ast} \varphi \mathd S_x, 
\end{align*}
where $\varphi$ is a test function, and $( \hat{n} \cdot F)^{\ast}$ is a
suitably chosen numerical flux in the unit normal direction $\hat{n}$.
Following \cite{hesthaven_nodal_2007}, we find a strong-DG form of
this system as
\begin{equation}
  0 = \int_{\D_k} u_t \varphi + [\nabla \cdot F (u)] \varphi \mathd x-
  \int_{\partial \D_k} [ \hat{n} \cdot F - ( \hat{n} \cdot F)^{\ast}] 
  \varphi \mathd S_x.
  \label{eq:strong-dg}
\end{equation}
We seek to find a numerical vector solution $u^k \assign u_N|_{\D_k}$ from the
space $P_N^n (\D_k)$ of local polynomials of maximum total degree $N$ on each
element.  We choose the scalar test function $\varphi \in P_N (\D_k)$ from the
same space and represent both by expansion in a basis of $N_p\assign
\operatorname{dim} P_N(\D_k)$ Lagrange polynomials $l_i$ with respect to a set
of interpolation nodes \cite{warburton_explicit_2006}. We define the mass,
stiffness, differentiation, and face mass matrices
\begin{subequations} 
  \label{eq:dg-global-matrices} 
  \begin{align} 
    M_{i j}^k & \assign \int_{\D_k} l_i l_j \mathd x, \\ 
    S_{i j}^{k, \partial\nu} & \assign \int_{\D_k} l_i \partial_{x_{\nu}} l_j \mathd x,\\ 
    D^{k, \partial\nu} & \assign (M^k)^{- 1} S^{k, \partial\nu},\\ 
    M_{i j}^{k, A} & \assign \int_{A \subset \partial \D_k} l_i l_j \mathd S_x.
  \end{align} 
\end{subequations}
Using these matrices, we rewrite \eqref{eq:strong-dg} as
\begin{align} 
  0 & =  M^k \partial_t u^k 
  + \sum_{\nu} S^{k, \partial_{\nu}} [F(u^k)]
  - \sum_{F \subset \partial \D_k} M^{k, A} [
  \hat{n} \cdot F - ( \hat{n} \cdot F)^{\ast}], \notag \\
  \label{eq:semidiscrete-dg} 
  \partial_t u^k & = 
  - \sum_{\nu} D^{k, \partial_{\nu}} [F(u^k)]
  + L^k [ \hat{n} \cdot F 
  - ( \hat{n} \cdot F)^{\ast}] |_{A \subset \partial \D_k}.
\end{align}
The matrix $L^k$ used in \eqref{eq:semidiscrete-dg} deserves a little
more explanation. It acts on vectors of the shape
$[u^k|_{A_1},\dots,u^k|_{A_4}]^T$, where $u^k|_{A_i}$ is the vector of facial
degrees of freedom on face $i$. For these vectors, $L^k$ combines the effect
of applying each face's mass matrix, embedding the resulting facial values
back into a volume vector, and applying the inverse volume mass matrix.
Since it ``lifts'' facial contributions to volume contributions, it is
called the \emph{lifting matrix}. Its construction is shown in Figure
\ref{fig:lifting-matrix}. 

\begin{myfigure}
\begin{center}
\begin{tikzpicture}[
  densemat/.style={fill=gray!30},
  ]
  \draw [densemat] (-1,-1) rectangle (1,1) ;
  \draw (2,-1) rectangle +(4,2) ;
  \draw [densemat] (-2,-1) rectangle (-5,1) ;

  \draw [densemat] (2,1) rectangle +(1,-1) ;
  \foreach \ys/\ye in {0.1/0.3, -0.25/-0.7, 0.6/0.8}
  { \draw [densemat] (3,\ys) rectangle (4,\ye) ; }
  \foreach \ys/\ye in {0.4/0.1, -0.1/-0.3, -0.9/-1, 1/0.9}
  { \draw [densemat] (4,\ys) rectangle (5,\ye) ; }
  \foreach \ys/\ye in {0.8/0.6, 0.3/0, -0.9/-0.7, -0.2/-0.5}
  { \draw [densemat] (5,\ys) rectangle (6,\ye) ; }

  \node at (2.5,0.5) {$M^{k,A_1}$} ;
  \node at (3.5,-0.475) {$M^{k,A_2}$} ;
  \node at (4.5,0.25) {$M^{k,A_3}$} ;
  \node at (5.5,-0.35) {$M^{k,A_4}$} ;

  \node at (0,0) {$(M^k)^{-1}$} ;
  \node at (1.5,0) {$\cdot$} ;
  \node at (-1.5,0) {$=$} ;
  \node at (-3.5,0) {$L^k$} ;

  \draw [|<->|] (6.25,-1) -- +(0,2) node [pos=0.5, anchor=west] {$N_p$} ;

  \draw [|<->|] (2,-1.25) -- (3,-1.25) node [anchor=west] {$N_{fp}$} ;
\end{tikzpicture}
\end{center}
\caption{Construction of the Lifting Matrix $L^k$.}
\label{fig:lifting-matrix}
\end{myfigure}

It deserves explicit mention at this point that the left multiplication by the
inverse of the mass matrix that yields the explicit semidiscrete scheme
\eqref{eq:semidiscrete-dg} is an elementwise operation and therefore feasible
without global communication. This strongly distinguishes DG from other finite
element methods. It enables the use of explicit (e.g., Runge-Kutta)
timestepping and greatly simplifies our efforts of bringing DG onto the GPU.

\subsection{Implementing DG}

\label{ssec:implement-dg}

DG decomposes very naturally into four stages, as visualized in Figure
\ref{fig:dg-subtasks}.  This clean decomposition of tasks stems from the fact
that the discrete DG operator \eqref{eq:semidiscrete-dg} has two additive
terms, one involving an element volume integral, the other an element surface
integral.  The surface integral term then decomposes further into a `gather'
stage that computes the term
\begin{equation} 
  [ \hat{n} \cdot F(u_N^-) - ( \hat{n} \cdot F)^{\ast}(u_N^-, u_N^+)]|_{A \subset
  \partial \D_k} 
  \label{eq:num-flux} 
\end{equation}
and a subsequent lifting stage. The notation $u_N^-$ indicates the value of 
$u_N$ on the face $A$ of element $\D_k$, $u_N^+$ the value of $u_N$ on the face 
opposite to $A$.

As is apparent from our use of a Lagrange
basis, we implement a \emph{nodal} version of DG, in which the stored degrees
of freedom (``\emph{DOF}s'') represent the values of $u_N$ at a set of
interpolation nodes. This representation allows us to find the facial values
used in \eqref{eq:num-flux} by picking the facial nodes from the volume
field. (This contrasts with a \emph{modal} implementation in which DOFs
represent expansion coefficients.  Finding the facial information to compute 
\eqref{eq:num-flux} requires a different approach in these schemes.)

Observe that most of DG's stages are \emph{element-local} in the sense that
they do not use information from neighboring elements. Moreover, these local
operations are often efficiently represented by a dense matrix-vector
multiplication on each element.

\begin{myfigure}
  \centering
  \begin{tikzpicture}[
    txt/.style={text height=1.5ex, text depth=0.25ex},
    operation/.style={rectangle,draw,minimum height=5ex,txt},
    localop/.style={operation,line width=2pt,txt},
    data/.style={circle,draw,minimum size=7ex,txt},
    ]

    \node [data] (state) { $u^k$ } ;
    \node [operation] (gather) [above right=0.1cm and 0.5cm of state] { Flux Gather } ;
    \node [localop] (lift) [right=1cm of gather] { Flux Lifting } ;
    \node [localop] (fu) [below right=0.1cm and 0.5cm of state] { $F(u^k)$ } ;
    \node [localop] (diff) [right=1cm of fu] { Local Differentiation} ;
    \node [data] (rhs) [above right=0.1cm and 0.5cm of diff] { $\partial_t u^k$ } ;
    \draw [->] (state) |- (gather) ;
    \draw [->] (gather) |- (lift) ;
    \draw [->] (lift) -| (rhs) ;

    \draw [->] (state) |- (fu) ;
    \draw [->] (fu) |- (diff) ;
    \draw [->] (diff) -| (rhs) ;
  \end{tikzpicture}
  \caption{Decomposition of a DG operator into
  subtasks. Element-local operations are highlighted with a bold outline.}
  \label{fig:dg-subtasks}
\end{myfigure}

It is worth noting that since simplicial elements only require affine
transformations $\Psi_k$ from reference to global element, the global matrices
can easily be expressed in terms of reference matrices that are the same for
each element, combined with scaling or linear combination, for example
\begin{subequations}
\label{eq:dg-reference-matrices}
\begin{align}
  M_{i j}^k & = 
  \underbrace{\left|\det \frac{\mathd \Psi_k}{\mathd r}\right|}_{J_k\assign}
  \underbrace{\int_{\mathsf{I}} l_i l_j \mathd x}_{M_{ij}\assign},\\
  S_{i j}^{k, \partial \nu} & = J_k \sum_{\mu} \frac{\partial \Psi_{\nu}}{\partial r_{\mu}}
  \underbrace{\int_{\mathsf{I}} l_i \partial_{r_{\mu}} l_j \mathd x}_{S_{ij}^{\partial \mu}\assign},
\end{align}
\end{subequations}
where $\mathsf{I}=\Psi_k^{-1}(\D_k)$ is a reference element. We define the remaining
reference matrices $D$, $M^A$, and $L$ in an analogous fashion.

\section{The CUDA Parallel Computation Model}

\label{sec:gpu-hardware}

Graphics hardware is aimed at the real-time rendering of large numbers of
textured geometric primitives, with varying amounts of per-pixel and
per-primitive processing. This problem is, for the most part, embarrassingly
parallel and exhibits this parallelism at both the pixel and the primitive
level.  It is therefore not surprising that the parallelism delivered by
graphics-derived computation hardware also exhibits two levels of parallelism.
On the Nvidia hardware \cite{lindholm_nvidia_2008} targeted in this work, up to
30 independent, parallel \emph{multiprocessors} form the higher level. Each of
these multiprocessors is capable of maintaining several hundred threads in
flight at any given time, giving rise to the lower level.

One such multiprocessor consists of eight functional units controlled by a
single instruction decode unit.  Each of the functional units, in turn, is
capable of executing one basic single-precision floating-point or integer
operation per clock cycle. Interestingly, a fused floating-point multiply-add
is one of these basic operations. The instruction decode unit feeding the eight
functional units is capable of issuing one instruction every four clock cycles,
and therefore the smallest scheduling unit on this hardware is what Nvidia
calls a \emph{warp}, a set of $T \assign 32$ threads. The architecture is
distinguished from conventional single-instruction-multiple-data (\emph{SIMD})
hardware by allowing threads within a warp to take different branches, although
in this case each branch is executed in sequence. To emphasize the difference,
Nvidia calls Tesla a single-instruction-multiple-thread (\emph{SIMT})
architecture.

Up to 16 of these warps are now aggregated into a \emph{thread block} and sent
to execute on a single multiprocessor. Threads in a block share a piece of
execution hardware, and are hence able to take advantage of additional
communication facilities present in a multiprocessor, namely, a barrier that
may optionally serve as a memory fence, and 16KiB\footnote{``KiB'' stands for
\emph{Kilobyte binary} or \emph{Kibibyte} and represents $1024=2^{10}$ bytes.
\cite{iec_kibi_2000}} of banked\footnote{``Banking'' of shared memory means
that only addresses in distinct banks can be accessed simultaneously. Allowing
simultaneous access to \emph{all} addresses in shared memory would require
prohibitive amounts of addressing logic. Therefore, banking is an expected
feature of parallel memory.} shared memory. The shared memory has 16 banks,
such that half a warp accesses shared memory simultaneously. If all 16 threads
access different banks, or if all 16 access the same memory location (a
\emph{broadcast}), the access proceeds at full speed. Otherwise, the whole warp
waits as maximal subsets of non-conflicting accesses are carried out
sequentially.

A potentially very large number of thread blocks is then aggregated into a
\emph{grid} and forms the unit in which the controlling host processor submits
work to the GPU. There is no guaranteed ordering between thread blocks in a
grid, and no communication is allowed between them.  Only after successful
completion of a grid submission, the work of all thread blocks is guaranteed to
be visible. In that sense, grid submission serves as a synchronization point.

Indices within a thread block and within a grid are available to the program at
run time and are permitted to be multi-dimensional to avoid expensive integer
divisions. We will refer to these indices by the symbols $t_x$, $t_y$, $t_z$,
and $b_x$, $b_y$.

All threads have read-write access to the GPU's on-board (`\emph{global}')
memory.  A single access to this off-chip memory has a latency of several
hundred clock cycles. To hide this latency, a multiprocessor will schedule
other warps if available and ready. A few things influence how many threads are
available: Each thread requires a number of registers. Also, the work of a
group of threads often involves a certain amount of shared memory. More threads
may therefore also consume more shared memory.  Since both the register file
and the amount of shared memory is finite, their use may lead to artificial
limits on the number of threads in a block. If there are very few threads in a
block and there isn't space for many blocks on the same multiprocessor, the
device may fail to find warps it can run while waiting for memory transactions.
This decreases global memory bandwidth utilization.  Another aspect influencing the
available bandwith to global memory is the pattern in which access occurs.
Taking 32-bit accesses as an example, loads and stores to global memory achieve
the highest bandwidth if, within a warp, thread $i$ accesses memory location $b
+ \pi(i)$, where $b$ is a 16-aligned base address and $\pi$ is a mapping
obeying $\lfloor\pi(i)/ 16\rfloor=\lfloor i / 16\rfloor$. Note that for global
\emph{fetches} only, these restrictions can be alleviated somewhat through the
use of \emph{texture units}.

A final bit of perspective: While the graphics card achieves an order of
magnitude larger bandwith to its global memory than a conventional processor
does to its main memory, its floating point capacity eclipses this already
large bandwidth by yet another order of magnitude. If we visualize both compute
and memory bandwidth as physical ``pipes'' with a certain diameter, the
challenge in designing algorithms for this architecture lies in keeping each
pipe flowing at capacity while using a minimum of buffer space.

\section{DG on the GPU: Design}

\label{sec:dg-gpu-design}

The answers to three questions emerge as the central
design decisions in mapping a numerical method into an algorithm that can run 
on a GPU:
\begin{description}
  \item [Computation Layout.] How can the task be decomposed into a grid of 
    thread blocks, given there cannot be any inter-block communication?
    Do we need a sequence of grids instead of a single grid?
  \item [Data Layout.] How well does the data conform to the device's alignment
    requirements? Where and to what extent will padding be used?
  \item [Fetch Schedule.] When will what piece of the data be fetched from 
    global into on-chip memory, i.e. registers or shared memory?
\end{description}

Note that the computation layout and the data layout are often the same, and
rarely independent. For the bandwidth reasons described in Section
\ref{sec:gpu-hardware}, the index of the thread computing a certain result
should match the index where that result is stored.  Post-computation
permutations come at the cost of setting aside shared memory to perform the
permutation. It is therefore common to see algorithms designed around the
principle of \emph{one thread per output}.  The fetch schedule, lastly,
determines how often data can be reused before it is evicted from on-chip
storage.

Unstructured discontinuous Galerkin methods have a number of natural
granularities: 
\begin{itemize}
  \item the number $N_p$ of DOFs per element,
  \item the number $N_{fp}$ of DOFs per face,
  \item the number $N_f$ of faces per element,
  \item the number $n$ of unknowns in the system of conservation laws.
\end{itemize}
The number of elements $K$ also influences the work partition, but it is less
important in the present discussion. 

The first three granularities above depend on the chosen order of approximation
as well as the shape of the reference element. Figure \ref{fig:dof-counts}
gives a few examples of their values. Perhaps the first problem that needs to
be addressed is that many of the DOF counts, especially at the practically
relevant orders of 3 and 4, conform quite poorly to the hardware's preference
for batches of 16 and 32.  A simple solution is to round the size of each
element up to the next alignment boundary. This leads to a large amount of
wasted memory.  More severely, it also leads to a large amount of wasted
processing power, assuming a one-thread-per-output design. For example,
rounding $N_p$ for a fourth-order element up to the next warp size boundary
($T=32$) leads to 45\% of the available processing power being wasted.  It is
thus natural to aggregate a number of elements to get closer to an alignment
boundary. Now, each of the parts of a DG operator is likely to have its
own preferred granularity corresponding to one thread block. One option is to
impose one such part's granularity on the whole method.  We find that a better
compromise is to introduce a sub-block granularity for this purpose. We
aggregate the smallest number $K_M$ of elements to achieve less than 5\% waste
when padding up to the next multiple $N_{pM}$ of $T/2=16$.  Figure
\ref{fig:microblock-mem-layout} illustrates the principle.  We then impose the
restriction that each thread block work on an integer number of these
\emph{microblocks}. We assign the symbol $n_M\assign \lceil K/K_M\rceil$ to the
total number of microblocks.

\begin{myfigure}
  \centering
  \subfigure[DOF counts for moderate-order tetrahedral elements.]{
    \label{fig:dof-counts}
    \begin{tabular}[b]{cccc}
    \hline
    $N$ & $N_p$ & $N_{fp}$ & $N_f N_{fp}$ \\
    \hline
    1 & 4 & 3 & 12 \\
    2 & 10 & 6 & 24 \\
    3 & 20 & 10 & 40 \\
    4 & 35 & 15 & 60 \\
    5 & 56 & 21 & 84 \\
    6 & 84 & 28 & 112 \\
    7 & 120 & 36 & 144 \\
    \hline
    \end{tabular}
  }
  \hspace{1em}
  \subfigure[Microblocked memory layout.]{
    \label{fig:microblock-mem-layout}
    \begin{tikzpicture}[scale=0.7]
      \foreach \bottom in {0, -1, -2}
      {
        \draw (0,\bottom) rectangle +(6.4,1) ;
        \draw [fill=gray!20] (0,\bottom) rectangle +(2,1) ;
        \draw [fill=gray!20] (2,\bottom) rectangle +(2,1) ;
        \draw [fill=gray!20] (4,\bottom) rectangle +(2,1) ;
      }
      \foreach \i in {0,0.1,...,6.4}
      { \draw (\i,1) -- (\i,0.9); }

      \foreach \x in {1,3,5}
      {
        \node at (\x,-1.5) {Element} ;
        \node at (\x,-0.5) {Element} ;
        \node at (\x,0.4) {\dots} ;
      }

      \node (padding) at (6,2.3) {Padding} ;
      \draw [->] (padding) -- (6.2,0.5) ;

      \draw [|<->|] (0,-2.5) -- +(2,0)
        node [pos=0.5,anchor=north] {$N_p$} ;
      \draw [|<->|] (0,-2.25) -- +(6,0) 
        node [pos=0.5,anchor=north] {$K_M N_p$} ;
      
      \node [left] at (-0.2,0.5) { 128 };
      \node [left] at (-0.2,-0.5) { 64 };
      \node [left] at (-0.2,-1.5) { 0 };
    \end{tikzpicture}
  }
  \caption{Matching DG granularities to GPU alignment boundaries.}
\end{myfigure}

The next question to be answered involves decomposing a task into an
appropriate set of thread blocks. This decomposition is problem-dependent, but
a few things can be said in general.  We assume a task that has to be performed
in parallel, independently, on a number of work units, and that it requires
some measure of preparation before actual work units can be processed. We are
trying to find the right amount of work to be done by a single thread block. We
may let the block complete work units in parallel, alongside each other in a
single thread (`\emph{inline}' for brevity), or sequentially. We will use the
symbols $w_p$, $w_i$ and $w_s$ for the number of work units processed in each
way by one thread block. Thus the total number of work units processed by one
thread block is $w_p w_i w_s$.  Increasing $w_p$ often through increased
parallelism and reuse of data in shared memory, but typically also requires
additional shared memory buffer space.  Increasing $w_i$ gains speed through
reuse of data in registers. For example, take a two-operand procedure like
matrix multiplication. Increasing $w_i$ allows a single thread to use data from
the first operand, once loaded into registers, to process more than one column
of the second operand. Like $w_p$, varying $w_i$ also influences buffer space
requirements.  $w_s$, finally, amortizes preparation work over a certain number
of work units, at the expense of making the computation more granular.
Achieving a balance between these aspects is not generally straightforward, as
Figure \ref{fig:diff4-mat-smem-work-distribution} will demonstrate. Note that
each of the methods discussed below will have its own values for $w_p$, $w_i$,
and $w_s$.

We noted above that the number $n$ of variables in the system of conservation
laws \eqref{eq:claw} also introduces a granularity. In some cases, it may be
advantageous to allow this system size to play a role in deciding data and
computation layouts. One might attempt do this by choosing a packed field
layout, i.e.  by storing all field values at one node in $n$ consecutive memory
locations.  However, a packed field layout is not desirable for a number of
reasons, the most significant of which is that it is unsuited to a
one-thread-per-output computation. If thread 0 computes the first field
component, thread 1 the second, and so on, then each field component is found
by evaluating a different expression, and hence by different code. This cannot
be efficiently implemented on SIMT hardware. One could also propose to take
advantage of the granularity $n$ by letting one thread compute all $n$
different expressions in the conservation law for one node. It is practical to
exploit this for the gathering of the fluxes and the evaluation of $F(u)$.  For
the more complicated lifting and differentiation stages on the other hand, this
leads to impractical amounts of register pressure. We find that, especially at
moderate orders, the extra flexibility afforded by ignoring $n$ outweighs any
advantage gained by heeding it. If desired, one can always choose $K_M=n$ or
$w_i=n$ to closely emulate the strategies above. Further, note that for the
linear case discussed here, one has significant freedom in the ordering of
operations, for example by commuting the evaluation of $F(u^k)$ with local
differentiation.

A final question in the overall algorithm design is whether it is appropriate
to split the DG operator into the subtasks indicated in Figure
\ref{fig:dg-subtasks}, rather than to use a single or only two grids to compute
the whole operator. Field data would need to be fetched only once, leading to a
good amount of data reuse. But at least for the scarce amounts of shared memory
buffer space in current-generation hardware, this view is too simplistic. Each
individual subtask tends to have a better, individual use for on-chip memory.
Also, it is tempting to combine the gather and lift stages, since one works on
the immediate output of the other.  Observe however that there is a mismatch in
output sizes between the two. For each element, the gather outputs $N_{fp}N_f$
values, while the lift outputs $N_p$. These two numbers differ, and therefore
the optimal computation layouts for both kernels also differ.  While it is
possible to use the larger of the two computation layouts and just idle the
overlap for the other computation, this is suboptimal. We find that the added
fetch cost is easily amortized by using an optimal computation layout for each
part of the flux treatment.

\section{DG on the GPU: Implementation}

\label{sec:dg-gpu-implementation}

\subsection{How to read this Section}
\label{sec:pseudocode-notation}

To facilitate a detailed, yet concise look at our implementation techniques,
this section supplements its discussion with pseudocode for some particularly
important subroutines. Pseudocode contains all the implementation details and
exposes the basic control and synchronization structure at a single glance. In
addition to the code, there is text discussing every important design decision
reflected in the code.

To maximize readability, we rely on a number of notational conventions. First,
$\intceil x n$ is the smallest integer larger than $x$ divisible by $n$. Next,
$\halfopen a b$ denotes the `half-open' set of integers $\{a,\ldots,b-1\}$.
Using this notation, we may indicate `vectorized' statements, e.g. an
assignment $a_{\halfopen k {k+n}} \leftarrow k_{\halfopen 0 n}$. The loops 
indicated by these statements are always fully unrolled in actual code. 
Depending on notational convenience,
we alternate between subscript notation $a_i$ and indexing notation
$a[i]$. Both are to be taken as equivalent. Sometimes, we use both sub- and
superscripts on a variable. This helps brevity and readability, but is only
done if the memory layout of the corresponding variable is clarified elsewhere.
Otherwise, for multidimensional indices, C-like (row-major) data layout is
assumed.

Lastly, the GPU offers many different types of storage. To avoid confusion, we
assign each type of storage a separate typographical convention, as outlined in
Table \ref{tab:storage-typo}. If and only if two storage locations of
different types are used for related data, we use the same letter for both.

\begin{table}
  \begin{center}
    \begin{tabular}{cll}
      \hline
      Convention & & Storage Type \\
      \hline
      $v$ & Italic font & Constant or unrolled loop variable \\
      $\reg{v}$ & Typewriter font & Register variable \\
      $\smem{v}$ & Superscript $\mathtt S$ & Variable in shared memory \\
      $\gmem{v}$ & Superscript $\mathtt G$ & Variable in global memory \\
      $\tex{v}$ & Superscript $\mathtt T$ & Variable bound to a texture \\
      \hline
    \end{tabular}
  \end{center}
  \caption{Typographical conventions for different types of GPU storage.}
  \label{tab:storage-typo}
\end{table}

\subsection{Flux Lifting}
\label{ssec:flux-lift}

Lifting is one of the \emph{element-local} components of a discontinuous
Galerkin operator, and, for simplicial elements, is efficiently represented by
a matrix-matrix multiplication as in Figure \ref{fig:el-local-matmul},
followed by an elementwise scaling. 

\begin{myfigure}
  \centering
  \subfigure[
    Applying an element-local DG operator $L$ to a field $u$ 
    by a matrix-matrix product.
  ]{
    \label{fig:el-local-matmul}
    \begin{tikzpicture}
      \newcommand{\elements}[1]{
        \draw [fill=gray!30] (0,0) rectangle (2.5,1.5) ; 
        \draw [xstep=0.2,ystep=2, gray] (0,0) grid (1.8,1.5) ; 
        \draw (0,0) rectangle (2.5,1.5) ;
        \node at (1.2,0.75) {#1} ;
        \node at (2.15,0.75) {\dots} ;
        \foreach \i in {0,1,2,3,8}
        { 
          \node 
          [rotate=270,font=\tiny,anchor=south, inner sep=0.05cm, text=gray] 
            at (0.2*\i,0.75) { Element \i }; 
        }
      }
      \node[matrix,ampersand replacement=\&, row sep=0.25cm, column sep=0.25cm, inner sep=0]
      {
        \&
        \elements{$u$}
        \\
        \draw [fill=gray!30] (0,0) rectangle (1.5,1.5) 
          node [pos=0.5] {$L$} ;
        \&
        \elements{$Lu$}
        \\
      };
    \end{tikzpicture}
  }
  \hspace{0.25em}
  \subfigure[
    Output memory layout for the flux gather stage, 
    input memory layout of the flux lifting stage.
  ]{
    \label{fig:gather-mem-layout}
    \begin{tikzpicture}[scale=0.45]
      \foreach \bottom in {0, -1, -2}
      {
        \draw (0,\bottom) rectangle +(12.8,1) ;
        \foreach \elleft in {0,4,8}
        { 
          \draw [fill=gray!20] (\elleft,\bottom) rectangle +(4,1) ; 
          \node [anchor=north] at ($(\elleft,\bottom)+(2,1.1)$) {\scriptsize Element } ;
          \foreach \face in {0,1,2,3}
          {
            \draw (\elleft+\face,\bottom) -- +(0,0.4) ;
            \node at ($(\elleft+\face,\bottom)+(0.5,0.2)$) {\tiny Face} ;
          }
        }
      }

      \foreach \i in {0,0.1,...,12.8}
      { \draw (\i,1) -- (\i,0.9); }

      \node (padding) at (12,2) {Padding} ;
      \draw [->] (padding) -- (12.4,0.5) ;

      \draw [|<->|] 
        (0,-3) -- +(1,0)
        node [pos=0.5,anchor=north] {$N_{fp}$} ;
      \draw [|<->|] 
        (0,-2.75) -- +(4,0)
        node [pos=0.55,anchor=north] {$N_f N_{fp}$} ;
      \draw [|<->|] 
        (0,-2.5) -- +(12,0)
        node [pos=0.5,anchor=north] {$K_M N_f N_{fp}$} ;
      \draw [|<->|] 
        (0,-2.25) -- +(12.8,0)
        node [pos=1,below left=0.1 and 0,inner xsep=0,outer xsep=0] {$N_{fM}$} ;
        
      \node [left] at (-0.2,0.5) { 256 };
      \node [left] at (-0.2,-0.5) { 128 };
      \node [left] at (-0.2,-1.5) { 0 };
    \end{tikzpicture}
  }
  \caption{Implementation aspects of flux lifting.}
\end{myfigure}

The first, tempting approach to implementing this is to take advantage of the
vendor-provided GPU-based BLAS workalike. This is hampered by sub-optimal
performance and strict alignment requirements. As a result, a custom algorithm
is in order.

One key to high performance on the GPU is to find a good use for the scarce
amount of shared memory. Both operands in an element-local matrix
multiplication see large amounts of reuse: Each field value is used $N_p$
times, and each entry of a local matrix is used $N_p$ times \emph{for each
element}. It is therefore a sensible wish to load both operands into shared
memory. For the tetrahedral elements targeted here, this is problematic.  Even
for elements of modest order, the matrix data quickly becomes too large. This
restricts the applicability of a matrix-in-shared approach to low orders, and
we will therefore first examine the more broadly applicable method of using the
shared memory for field data. Still, matrix-in-shared does provide a benefit
for certain low orders and is examined in the context of element-local
differentiation in Section \ref{ssec:local-diff}.

We choose a one-thread-per-output design for flux lifting. This dictates that
computation and output layouts match Figure \ref{fig:microblock-mem-layout}.
But the input layout for lifting is mildly different: The flux gather, which
provides the input to lifting, extracts $N_f N_{fp}$ DOFs per element. Recall
that the layout of Figure \ref{fig:microblock-mem-layout} provides $N_p$ DOFs
per element. Since typically $N_p \ne N_f N_{fp}$, we introduce a mildly
different layout as shown in Figure \ref{fig:gather-mem-layout}, using the same
number $K_M$ of elements as found in a mircroblock, padded to half-warp
granularity.  This padding is likely somewhat more wasteful than the carefully
tuned one of Figure \ref{fig:microblock-mem-layout}.  Fortunately, this is
irrelevant: We will not be using Figure \ref{fig:gather-mem-layout} as a
computation layout, and data in this format is used only for short-lived
intermediate results.  Overall, the resulting memory layout has
$N_{fM}\assign\intceil{N_f N_{fp} K_M}{T/2}$ DOFs per microblock.

We are now ready to discuss the actual algorithm, at the start of which we need
to fetch field data into shared memory. Because we chose a
one-thread-per-output computation layout, we will have $N_p$ threads per
element fetching data. Due to the mismatch between $N_p$ and $N_fN_{fp}$, we
may require multiple fetch cycles to fetch all data. In addition, the last such
fetch cycle must involve a length check to avoid overfetching. It is important
to unroll this fetch loop and to use some care with the ending conditional to
still allow fetch pipelining\footnote{Pipelining is a fetch optimization
strategy. It performs high-latency fetches in batches ahead of a computation.
Since a warp only stalls when unavailable data is actually used in a
computation, this allows a single thread to wait for multiple memory
transactions simultaneously, decreasing latency and reducing the need for
parallel occucpancy. The Nvidia compiler automatically pipelines fetches if the
code structure allows it.} to occur.

With the field data in shared memory, the matrix data is fetched using texture
units.  By way of the texture cache, we hope to take advantage of the
significant redundancy in these fetches. The matrix texture should use
column-major order: Realize that within a block, a large number of threads,
each assigned to a row of the matrix, load values from each column in turn.
Column-major order gives the most locality to this access pattern. 

With this preparation, the actual matrix-matrix product can be performed.
Since all threads within one element load each of the element's nodal values
from shared memory in order, these accesses are handled as a broadcast and
therefore conflict-free.  Conflicts do occur, however, if a half-warp straddles
an element boundary within a microblock. In that case, threads before and after
the element boundary access different elements, and therefore a
double-broadcast bank conflict occurs. Figure
\ref{fig:smem-field-comp-layout-bad} shows the genesis of this conflict.
Fortunately, that does not automatically mean that microblocking is a bad idea.
It turns out that the performance lost when using no microblocking and hence
full padding is about the same as the one lost to these bank conflicts.  Even
better: there is a way of mitigating the conflicts' impact \emph{without}
having to forgo the performance benefits of microblocking. The key realization
is that even if only one half of a warp encounters a conflict, the other half
of the warp is made to wait, too, regardless of whether it conflicted.
Conversely, if we assemble warps in such a way that conflict-prone and
non-conflict-prone half-warps are kept separate, then we avoid unnecessary
stalling. If $w_p>1$, then we can achieve such a grouping by laying out the
computation as seen in Figure \ref{fig:smem-field-comp-layout}.

\begin{myfigure}
  \centering
  \newcommand{\accessmublockback}
  {
    \draw [fill=gray!30] (mublock) rectangle (4,1) ;
    \draw (mublock) rectangle (4.8,1) ;
    \draw [xstep=2,ystep=0.5] (mublock) grid (4.8,1) ;
  }
  \newcommand{\accessmublockfront}
  {
    \foreach\halfwarp in {1.6,3.2}
    { 
      \draw [thick,densely dotted,red]
        (mublock) ++(\halfwarp,0) -- ++(0,1) ; 
    }

    \draw [|<->|] (mublock) ++(-0.25,0) -- ++(0,1) 
      node [pos=0.5, anchor=east] {$w_p$};
    \draw [|<->|] (mublock) ++(0,-0.25) -- ++(4.8,0) 
      node [pos=0.95, anchor=north] {$N_{pM}$};
    \draw [|<->|] (mublock) ++(0,-0.5) -- ++(4,0) 
      node [pos=0.75, anchor=north] {$N_pK_M$};
    \draw [|<->|] (mublock) ++(0,-0.75) -- ++(2,0) 
      node [pos=0.5, anchor=north] {$N_p$};


    \coordinate (threads) at ($ (mublock) + (0.05,0.3333) $);
    \foreach \p in {0,1}
    {
      \foreach \thread in {0,...,39}
      { \draw [fill=black] (threads) ++(0.1*\thread,0.5*\p) circle (0.015) ; }
    }

    \coordinate (eldata) at (1,2) ;
    \node [left=0.25cm,font=\footnotesize] 
      at ($ (eldata) + (0,0.5) $) {El. Data} ;

    \draw (eldata) rectangle +(2.2,1) ;
    \draw [xstep=1,ystep=0.5] (eldata) grid +(2.2,1) ;

    \foreach \p in {1}
    {
      \foreach \thread in {16,...,31}
      {
        \draw [->,red, thick, shorten <=0.05cm]
          (threads) ++(\thread*0.1,\p*0.5)
          .. controls +(0,0.5) and +(0,-0.5) .. 
          ($ (eldata) + ({0.25+floor(\thread/20)}, 0.25+0.5*\p) $)
          ;
      }
    }


    \coordinate (threadnum) at ($ (threads) +(0,-0.1666) $) ;
  }
  \tikzset{
    noconflict/.style={fill=darkgreen,opacity=0.2},
    hasconflict/.style={fill=red,opacity=0.6},
    }
  \subfigure[
    `Conventional', conflict-aggravating layout.
    The first and third warp (red) serialize access because of conflicts
    in the second half-warp of each microblock. Only the second warp (green)
    proceeds without conflicts.
  ]{
    \label{fig:smem-field-comp-layout-bad}
    \begin{tikzpicture}[scale=0.7]
      \coordinate (mublock) at (0,0) ;

      \accessmublockback

      \draw [hasconflict]
        (mublock) rectangle +(3.2,0.5) 
        ++(1.6,0.5) rectangle +(3.2,0.5) 
        ;

      \draw [noconflict]
        (mublock) ++(3.2,0) rectangle +(1.6,0.5) 
        (mublock) ++(0,0.5) rectangle +(1.6,0.5) 
        ;

      \accessmublockfront

      \draw [->,thick]
        (threadnum) -- +(47*0.1,0) 
        ++(0,0.5)-- +(47*0.1,0) 
        node [above=0.2cm,font=\footnotesize,text width=1cm] {thread number} ;
      \draw [thin]
        (threadnum) ++(47*0.1,0) 
        .. controls +(0,0.4) and +(0,-0.4) ..
        ++(-47*0.1,0.5) ;
    \end{tikzpicture}
  }
  \hspace{0.25em}
  \subfigure[
    Improved, conflict-mitigating layout.
    Only the second warp (red) serializes access for conflicts.
    The first and third warp (green) remain conflict-free.
  ]{
    \label{fig:smem-field-comp-layout}
    \begin{tikzpicture}[scale=0.7]
      \coordinate (mublock) at (0,0) ;

      \accessmublockback

      \draw [noconflict]
        (mublock) rectangle +(1.6,1) 
        ++(3.2,0) rectangle +(1.6,1) 
        ;

      \draw [hasconflict]
        (mublock) ++(1.6,0) rectangle +(1.6,1) 
        ;

      \accessmublockfront

      \foreach \hw in {0,1,2}
      {
        \draw [thick]
          (threadnum) 
          ++(\hw*16*0.1,0) -- +(15*0.1,0) 
          ++(0,0.5)-- +(15*0.1,0) 
          ;
        \draw [thin]
          (threadnum)
          ++({(\hw*16+15)*0.1},0) 
          .. controls +(0,0.4) and +(0,-0.4) ..
          ++(-15*0.1,0.5) ;
      }

      \foreach \hw in {0,1}
      {
        \draw [thin]
          (threadnum) 
          ++({(\hw*16+15)*0.1},0.5) 
          .. controls +(0,-0.4) and +(0,0.4) ..
          ++(1*0.1,-0.5) ;
      }

      \draw [thick,->]
        (threadnum)
        ++(3.2,0.5) -- ++(15*0.1,0)
        node [above=0.2cm,font=\footnotesize,text width=1cm] 
        {thread number} ;
    \end{tikzpicture}
  }
  \caption{Computation layouts for matrix multiplication
    with fields in shared memory.}
\end{myfigure}

Algorithm \ref{alg:lift} represents the aggregate of the techniques described
in this section. Observe that since there is no preparation work, we set
$w_s\assign 1$.  We should stress at this point that both the field-in-shared
and the matrix-in-shared approach can be used for both lifting and
element-local differentiation. Adapting the strategy of Algorithm
\ref{alg:lift} for the latter is quite straightforward.

\begin{algorithm}
\caption{Flux Lifting, field-in-shared.}
\label{alg:lift}
\begin{algorithmic}
  \REQUIRE A grid of $\intceil{n_M / w_p w_i }{} \times 1$
    blocks of size $T/2 \times w_p \times N_{pM} / (T/2)$.
  \REQUIRE Inputs:
    $\tex L$, the reference element's lifting matrix;
    $\tex i$, the per-element inverse Jacobians;
    $\gmem f$, the surface fluxes in the format of Figure 
    \ref{fig:gather-mem-layout}.
  \ENSURE Output: $\gmem r$, the surface fluxes 
    $\gmem f$ multiplied by the
    per-element lifting matrix $L^k$.
  \STATE $\reg m \leftarrow (b_x w_p+t_y)w_i$ 
    \COMMENT{the base microblock number}
  \STATE $\reg i \leftarrow (T/2) t_z + t_x$ 
    \COMMENT{this thread's DOF number within its microblock}
  \STATE \COMMENT{load data}
  \FORALL{unrolled 
  $b\in\halfopen{0}{\intceil{\intceil{N_{fM}}{T}/\intceil{N_{pM}}{T}}{}}$}
    \IF{$bN_{pM}+\reg i<N_{fM}$}
      \STATE $\smem f_{t_y,\halfopen{0}{w_i},bN_{pM} + \reg i} 
            \leftarrow \gmem f_{(\reg m+\halfopen{0}{w_i} )
            \intceil{N_{fM}}{T} +b N_{pM} + \reg i}$
    \ENDIF
  \ENDFOR
  \BARRIER
  \STATE \COMMENT{perform matrix multiply}
  \IF{$\reg i < K_M N_{pM}$}
    \STATE $\reg r_{\halfopen0{w_i}} \leftarrow 0$
    \FORALL{unrolled $n\in\halfopen 0 {N_f N_{fp}}$} 
      \STATE $\reg r_{\halfopen 0 {w_i}} 
      \leftarrow \reg r_{\halfopen 0 {w_i}} 
      + \tex L[\reg i \bmod N_p,n]\smem f_{t_y,\halfopen 0{w_i},n}$
    \ENDFOR
    \STATE $\gmem r_{(\reg m+\halfopen 0{w_i})N_{pM}+\reg i}
      \leftarrow 
      \tex i[(\reg m+\halfopen0{w_i})K_M+\lfloor \reg i / N_p\rfloor]
      \reg r_{\halfopen 0{w_i}}$
  \ENDIF
\end{algorithmic}
\end{algorithm}

\subsection{Flux Extraction}
\label{ssec:flux-gather}

In a strong-form, nodal implementation of the discontinuous Galerkin method,
flux extraction or `gather' iterates over the node indices of each face in the
mesh and evaluates the flux expression \eqref{eq:num-flux} at each such node.
As such, it is a rather quick operation characterized by few arithemtic
operations and a very scattered fetch pattern.  This non-local memory access
pattern is the most expensive aspect of flux extraction on a GPU, and our
foremost goal should therefore be to minimize the number of fetches at all
costs. For linear conservation laws, we may with very little harm treat the
element-local parts of a DG operator as if they acted on scalar fields.  This
is however not true of the non-local flux extraction. Fetching all fields only
once and then computing all $n$ fluxes saves a significant $n^2-n$ fetches of
each facial node value.

The next potential savings comes from the fact that the fluxes on the two sides of
an interior face pair use the same face data. By computing fluxes for such face
pairs together, we can cut the number of interior face fetches in half.  
Computing and storing opposite fluxes together is of course only 
possible if the task decomposition assigns both to the same thread block. We
will therefore need to invest some care into this decomposition.

To help find the properties of the task decomposition, observe that by choosing
to compute opposite fluxes together, we are implicitly rejecting a
one-thread-per-output design.  To accommodate opposite faces' fluxes being
computed simultaneously, we will allow the gathered fluxes to be written into a
shared memory buffer in random order in time, but conforming to the output
layout of Figure \ref{fig:gather-mem-layout}.  Once completed, this shared
memory buffer will then be flushed to global memory in one contiguous write
operation.  This limits our task decomposition choices: Thread blocks will
output contiguous pieces of data in the output layout.  This means that the
smallest granularity on which a thread block for flux extraction may begin and
end is that of a microblock: We will let each thread block compute fluxes on an
integer number $M_B$ of microblocks. Observe that this is not ideal: The
natural task decomposition for flux extraction is by face pair, not by element,
nor, even worse, by a group of elements as large as a microblock. Nonetheless,
given our output memory layout, this decomposition is inevitable.

But all is not lost. By carefully controlling the assignment of elements to
microblocks, and again by carefully choosing the assignment of microblocks to
flux extraction thread blocks, we can hope to recover many block-interior face
pairs within a thread block. Note the far-reaching consequences of what was
just decided: We need to have the elements participating in a flux-gather
thread block sit adjacent to each other in the mesh. To achieve this, we
partition the mesh into pieces of at most $K_M M_B$ elements each and then
assign the elements in each piece to microblocks sequentially. This means
nothing less than letting our mesh numbering be decided by what is convenient
for the gathering of fluxes.

What can we say about the required partition? It is important to realize that
this is a fairly different domain decomposition problem than the one for
distributed-memory machines.  First, there is a hard limit of $K_M M_B$
elements per piece, as determined by the amount of shared memory set aside for
write buffering.  Second, there is a (somewhat softer) limit on the number of
block-external faces. This limit stems from the fact that information about the
faces on which we gather fluxes needs to be stored somewhere.  Obviously,
block-internal face pairs can share this information and therefore require less
storage--one descriptor for each two faces. Face pairs on a block boundary are
less efficient. They require one descriptor for each face. If the block size
$K_M M_B$ is relatively large, a bad, splintered partition may have too many
boundary faces and therefore exceed the ``soft'' limit on available space for
face pair descriptors. Therefore, for large blocks, we require a `good'
partition with as few block-exterior face pairs as possible.  For very small
blocks, on the other hand, the problem is exactly opposite: If $K_M M_B$ is
small, the absolute quality of the mesh partition is not as critically
important: The small overall number of faces means that we will not run out of
descriptor space, making the soft limit even softer. 

So, how can the needed partition be obtained?  A natural first idea is to use
conventional graph partitioning software (e.g.  \cite{karypis_fast_1999}).
Problematically, these packages tend to fail when partitioning very large
meshes into very many small parts. In addition, our `soft' and `hard' limits
are difficult to enforce in these packages, so that obtaining a conforming
partition may take several `attempts' with increasing target partition sizes.
Increased target partition sizes, in turn, mean that there are microblocks
where element slots go unassigned.  This means that generic graph partitioners
are not a universal answer. They work well and generate good-quality partitions
if $K_M M_B \gtrapprox 10$. Otherwise, we fall back on a simple greedy
breadth-first agglomerator that picks elements by a total connectivity
heuristic, illustrated in Algorithm \ref{alg:greedy-partition}. In this case,
the greedy algorithm may produce a few very `bad' scattered blocks with many
external faces, but we have found that they matter neither in performance, nor
in keeping the `soft' limit.

\begin{algorithm}
\begin{algorithmic}
  \REQUIRE Input: set of elements $E$ with connectivity $C\assign\{ 
    (e_1,e_2):\text{$e_1$ and $e_2$ share a face}\}$.
  \ENSURE Output: the partition, a set of blocks $P$, each of size $\le l$.
    \STATE $P \leftarrow \emptyset$
    \WHILE{$E\ne\emptyset$}
      \STATE $Q \leftarrow \{\text{a seed element from $E$}\}$
        (a queue of candidate elements)
      \STATE $B \leftarrow \emptyset$
        (the block currently being generated)
      \LOOP
        \STATE Find and remove the element 
          $e\in Q$ that shares the most faces with $B$.
        \IF{$e\in E$}
          \STATE Remove $e$ from $E$, add it to $B$.
          \IF{$|B|=l$}
            \STATE Make first entry of $Q$ the new seed element, break the loop.
          \ENDIF
          \STATE $Q \leftarrow Q \cup \{f:(e,f)\in C\}$
        \ENDIF
        \IF{$Q=\emptyset$}
          \IF{$E=\emptyset$}
            \STATE Break the loop.
          \ELSE
            \STATE Add an arbitrary element from $E$ to $Q$.
          \ENDIF
        \ENDIF
      \ENDLOOP
      \STATE $P \leftarrow P \cup \{B\}$
    \ENDWHILE
\end{algorithmic}
\caption{Simple Greedy Partition.}
\label{alg:greedy-partition}
\end{algorithm}

Once the partition is constructed, we obtain for each block a number of
elements whose faces fall into one of three categories: intra-block interior,
inter-block interior, and boundary faces. We design our algorithm to walk an
array of data structures describing face pairs, each of which falls into one of
these categories. Within this array, each face pair structure contains all
information needed to gather and compute the fluxes for its target face(s).
Descriptors for intra-block interior face pairs drive the flux computation for
two faces at once, while the other two kinds only drive the computation for one
face.  The array is loaded from global into shared memory when each thread
block begins its work. To minimize branching and to save storage space in each
descriptor, we make the kind of each face pair descriptor implicit in its
position in the array. To achieve this, we order the array by the face pair's
category and store how many face pairs of each category are contained in the
array.

Because we implement a nodal DG method, face index lists play an important role
in the gather process: Each face's nodal values need to be extracted from a
given volume field. Since a tetrahedron has four faces, there are four possible
index subsets at which each face's DOFs are found, all of length $N_{fp}$.
Knowing these index subsets enables us to find surface nodal values for one
element. But we need to find \emph{corresponding} nodal values on two opposite
elements.  Therefore, we may need to permute the fetch ordering of one of the
elements in a face pair.  Altogether, to find opposing surface nodal values, we
need to store two index lists.  Since the number of distinct index lists is
finite, it is reasonable to remove each individual index list from the face
pair data structure and to instead refer to a global list of index lists.  We
find that a small texture provides a suitable storage location for this list.
Finally, note that intra-block face pairs require another index list: If we
strive to conform to an assumed `natural' face ordering of one `dominant' face,
writing the other's data into the purely facial structure from Figure
\ref{fig:gather-mem-layout} requires a different index list than the one needed
to read the element's volume data.

Of all the parts of a DG operator, the flux gather stage is the one that is
perhaps least suited to execution on a GPU. The algorithm is data-driven and
therefore branch-intensive, it accesses memory in an erratic way, and, as $n$
grows, it tends to require a fair bit of register space. It is encouraging to
see that despite these issues, it is possible to design a method, given in
Algorithm \ref{alg:flux-extraction}, that performs respectably on current
hardware. 

\begin{algorithm}
\caption{Flux Extraction.}
\label{alg:flux-extraction}
\begin{algorithmic}
  \newcommand{\lookup}[1]{\mathtt{.#1}}
  \REQUIRE a grid of $\lceil n_M/M_B\rceil \times 1$
    blocks of size $N_{fp} \times w_p\times 1$.
  \REQUIRE Inputs:
    $(\tex u)^{\halfopen 0 n}$, the set of fields of which fluxes are to be
    computed, each as a separate texture,
    $\gmem d$, face information records,
    $\tex J$, face index list array.
  \ENSURE Outputs:
    $(\gmem f)^{\halfopen 0 n}$, the surface fluxes for each face of each element,
    as a sequence of scalar fields.
  \STATE Load face information records from $\gmem d[b_x]$ into the shared memory variable $\smem d$.
  \BARRIER
  \STATE $\reg e\leftarrow t_y$ 
    \COMMENT{initialize the number of the face pair this thread is working on }
  \WHILE{$\reg e < \text{\# of interior face pairs in $\smem d$}$}
    \STATE $(\reg i^-, \reg i^+) 
      \leftarrow 
      \smem d[\reg e]\lookup{fetch\_base}^{-,+}
      +\tex J[\smem d[\reg e]\lookup{fetch\_idx\_list\_nr}^{-,+},t_x]$
    \STATE $
       \reg{u}^{\halfopen 0 n}_{-,+}
       \leftarrow (\tex u)^{\halfopen 0 n}_{\reg i^{-,+}}$
    \STATE ($\smem f)^{\halfopen 0 n}[\smem d[\reg e]\lookup{store\_base}^-+t_x]$\\
      $\qquad\leftarrow\smem d[\reg e]\lookup{face\_jacobian}\cdot
        [ \hat{n} \cdot F - ( \hat{n} \cdot F)^{\ast}]^{\halfopen 0 n}
        (\reg u^{\halfopen 0 n}_-, \reg u^{\halfopen 0 n}_+) $
    \STATE ($\smem f)^{\halfopen 0 n}[\smem d[\reg e]\lookup{store\_base}^+
      +\tex j[\smem d[\reg e]\lookup{store\_idx\_list\_nr}^+, t_x]]$\\
      $\qquad\leftarrow\smem d[\reg e]\lookup{face\_jacobian}\cdot
        [ (-\hat{n}) \cdot F - ( (-\hat{n}) \cdot F)^{\ast}]
        (\reg u^{\halfopen 0 n}_+, \reg u^{\halfopen 0 n}_-) $
    \STATE $\reg e \leftarrow \reg e + w_p$
  \ENDWHILE
  \WHILE{$\reg e < \text{\# of interior and exterior face pairs in $\smem d$}$}
    \STATE $(\reg i^-, \reg i^+) 
      \leftarrow 
      \smem d[\reg e]\lookup{fetch\_base}^{-,+}
      +\tex J[\smem d[\reg e]\lookup{fetch\_idx\_list\_nr}^{-,+},t_x]$
    \STATE $
       \reg{u}^{\halfopen 0 n}_{-,+}
       \leftarrow (\tex u)^{\halfopen 0 n}_{\reg i^{-,+}}$
    \STATE $(\smem f)^{\halfopen 0 n}[\smem d[\reg e]\lookup{store\_base}^-+t_x]$\\
      $\qquad\leftarrow\smem d[\reg e]\lookup{face\_jacobian}\cdot
        [ \hat{n} \cdot F - ( \hat{n} \cdot F)^{\ast}]
        (\reg u^{\halfopen 0 n}_-, \reg u^{\halfopen 0 n}_+) $
    \STATE $\reg e \leftarrow \reg e + w_p$
  \ENDWHILE
  \WHILE{$\reg e < \text{\# of face pairs in $\smem d$}$}
    \STATE $\reg i^-
      \leftarrow 
      \smem d[\reg e]\lookup{fetch\_base}^{-}
      +\tex J[\smem d[\reg e]\lookup{fetch\_idx\_list\_nr}^{-},t_x]$
    \STATE $
       \reg{u}^{\halfopen 0 n}_{-}
       \leftarrow (\tex u)^{\halfopen 0 n}_{\reg i^{-}}$
    \STATE $
       \reg{u}^{\halfopen 0 n}_{+}
       \leftarrow b(\reg u^{\halfopen 0 n}_-, \smem d[\reg e])$
       \COMMENT{calculate boundary condition}
    \STATE $(\smem f)^{\halfopen 0 n}[\smem d[\reg e]\lookup{store\_base}^-+t_x]$\\
      $\qquad\leftarrow\smem d[\reg e]\lookup{face\_jacobian}\cdot
        [ \hat{n} \cdot F - ( \hat{n} \cdot F)^{\ast}]
        (\reg u^{\halfopen 0 n}_-, \reg u^{\halfopen 0 n}_+) $
    \STATE $\reg e \leftarrow \reg e + w_p$
  \ENDWHILE
  \BARRIER
  \STATE 
    $(\gmem f)^{\halfopen 0 n}_{b_x M_B N_{fM} + \halfopen{0}{M_B N_{fM}}} \leftarrow
    (\smem f)^{\halfopen 0 n}_{\halfopen{0}{M_B N_{fM}}}$ (not unrolled)
\end{algorithmic}
\end{algorithm}

\subsection{Element-Local Differentiation}
\label{ssec:local-diff}

Unlike lifting, element-local differentiation must be represented not as
one matrix-matrix product (see Figure \ref{fig:el-local-matmul}), but as $d=3$
separate ones whose results are linearly combined to find the global $x$-, $y$-
and $z$-derivatives. Each of the $d$ differentiation matrices has $N_p\times N_p$
entries and is applied to the same data. To maximize data reuse and minimize
fetch traffic, it is immediately apparent that all $d$ matrix
multiplications should be carried out ``inline'' along with each other.

Superficially, this makes differentiation look quite like a lift where we have
chosen $w_i=d$. But there is one crucial difference: the three matrices used for
differentiation are all different. Increasing $w_i$ drives data reuse in
lifting simply by occupying more registers. As we will see in Section
\ref{sec:experiments}, this suffices to make it go very fast.  Differentiation
on the other hand already has a built-in ``$w_i$ multiplier'' of $d$ \emph{and}
has to deal with different matrices. Both factors significantly increase
register pressure.  Stated differently, this means that it is unlikely that we
will be able to drive matrix data reuse by using more registers as we were able
to do for lifting.  But the matrix remains the most-reused bit of data in the
algorithm.  In this section, we will therefore attempt to exploit this reuse by
storing the matrix, not the field, in shared memory.

We have already discussed in Section \ref{ssec:flux-lift} that the
matrix-in-shared approach can only work for low orders because of the rapid
growth of the matrix data with $N$. At first, this seems like a problematic
restriction that makes the approach less general than it could be. It can
however be turned into an advantage: Since we can assume that the algorithm
runs at orders six and below, we can exploit this fact in our design decisions.

We begin our discussion of this approach by figuring how the matrix data should
be loaded into shared memory. As in Section \ref{ssec:flux-lift}, we adopt a
one-thread-per-output approach. A straightforward first attempt may be to load
all $d$ local differentiation matrices into shared memory in their entirety.
Then each thread computes a different row of the matrix-vector product, and in
doing so, thread number $i$ accesses the $i$th row of the matrix.  Without loss
of generality, let the matrix be stored in row-major order, so that thread $i$
accesses memory cell number $i N_p$.  Shared memory has $T/2=16$ distinct
memory banks, and therefore the access is conflict-free iff $N_p$ and 16 are
relatively prime, or, more simply, iff $N_p$ is odd. This is encouraging: We
can achieve a conflict-free access pattern simply by adding a `padding' column
if necessary to enforce an odd stride $S$. Figure
\ref{fig:smem-matrix-assignment} shows the resulting assignment of matrix data
to shared memory banks, and Figure \ref{fig:smem-matrix-conflict-free}
illustrates the resulting conflict-free access pattern.

Unfortunately, this is too easy. In the presence of microblocking,
conflict-free access becomes more difficult.  If a half-warp straddles one or
more element boundaries, bank conflicts are likely to result.  The access
not only has a stride $S$, but also incorporates a jump from the end of the
matrix to its beginning, a stride of $(N_p-1) S$. And unlike in the
previous case, we cannot simply add a pad row to make the access 
conflict-free.  Figure \ref{fig:smem-matrix-conflicting} displays the problem.  

One way to avoid the disastrous end-to-beginning jump and to maintain the
conflict-free access pattern would be to duplicate the matrix data from the
first rows beyond the end of the matrix. This is workable in principle, but in
practice we are already filling the entire shared memory space with matrix data
and are unlikely to be able to afford the added duplication.  Fortunately, the
duplication idea can be saved, and there exists a conflict-free matrix storage
layout that does not require us to abandon microblocking.

Departing from the idea that we will store the \emph{entire} matrix, we aim at
storing just a constant-size row-wise \emph{segment} of the matrix.  Then, if
the end of the matrix falls within a segment, we fill up the rest of the
segment with rows from the beginning, providing the necessary duplication for
conflict-free access.  For this layout, we consider a composite matrix made up
of $N_M$ vertically concatenated copies of the $D^{\partial\mu}$. This
composite matrix is then segmented into pieces of $N_R$ rows each, where $N_R$
is chosen as a multiple of $T/2$.  Each such matrix segment has a naturally
corresponding range of degrees of freedom in a microblock, and we limit the
thread block that loads this matrix segment to computing outputs from this
range. Figure \ref{fig:smem-matrix-segment} illustrates the principle.

\begin{myfigure}
  \begin{center}
    \begin{tikzpicture}[scale=0.7]

      \coordinate (fld) at (3.5,3.3+3.2) ;
      \coordinate (res) at (3.5,3) ;
      \coordinate (mtx) at (0,3) ;

      \foreach \chunk in {0,0.8,...,3}
      { 
        \draw [red,fill=red!20] (0,3-\chunk) rectangle +(3.2,-0.8); 
        \draw [red,fill=red!20] ($(res) - (0,\chunk)$) rectangle +(3,-0.8); 
      }

      \draw (mtx) rectangle +(3.2,-3.2) ;

      \foreach \k in {0,1,2} 
      {
        \draw [fill=gray!20,opacity=0.7] ($(0,3)+\k*(1,-1)$) rectangle +(1,-1) ;
        \draw [fill=gray!20] ($(fld)-(0,\k)$) rectangle +(3,-1) ;
        \draw [fill=gray!20,opacity=0.7] ($(res)-(0,\k)$) rectangle +(3,-1) ;
      }

      \draw [fill=white] ($(fld) + (0,-3)$) rectangle +(3,-0.2) ;
      \draw [dashed] (0,0) -- (3,0) -- (3,3) ;
      \foreach \k in {0,0.2,...,3}
      {
        \draw ($(fld) + (\k,0)$) -- +(0,-3.2);
        \draw ($(res) + (\k,0)$) -- +(0,-3.2);
      }
      \foreach \k in {0,1,2} 
      { \node at ($(0.5,2.5) +\k*(1,-1)$) {$D$} ; }

      \draw [|<->|] 
        ($(mtx) + (-0.2,0)$) 
        -- ++(0,-0.4) node [anchor=east] {$N_R$} 
        -- ++(0,-0.4);

      \draw [|<->|] 
        ($(mtx) + (0,0.2)$) 
        -- ++(0.5,0) node [anchor=south] {$N_p$} 
        -- ++(0.5,0);

      \draw [|<->|] 
        ($(res) + (3.2,0)$) 
        -- ++(0,-3.2) node [anchor=west,pos=1.05] {$N_{pM}$} ;

      \draw [|<->|] 
        ($(res) + (3.4,0)$) 
        -- ++(0,-3) node [anchor=west,pos=0.5] {$N_M N_p$} ;

      \draw [|<->|] 
        ($(res) + (3.6,0)$) 
        -- ++(0,-1) node [anchor=west,pos=0.5] {$N_p$} ;

      \node at ($(fld) + 0.5*(3.2,-3)$) {\large $u$} ;
      \node at ($(res) + 0.5*(3.2,-3)$) {\large $Du$} ;

      \draw [opacity=0.3,fill=blue] 
        (res) rectangle +(0.2,-0.8) 
        coordinate [pos=0.5] (res-seg-1);

      \draw [opacity=0.3,fill=blue] 
        (res) ++(0.2,-0.8) rectangle +(0.2,-0.8);

      \draw [opacity=0.3,fill=darkgreen] 
        (fld) rectangle +(0.2,-1) 
        coordinate [pos=0.5] (fld-el-1-1);

      \draw [opacity=0.3,fill=darkgreen] 
        (fld) ++(0.2,0) rectangle +(0.2,-1) 
        coordinate [pos=0.5] (fld-el-2-1);
      \draw [opacity=0.3,fill=darkgreen] 
        (fld) ++(0.2,-1) rectangle +(0.2,-1) 
        coordinate [pos=0.5] (fld-el-2-2);

      \draw [thick,blue,->]
        (res-seg-1)
        .. controls +(-1,0) and +(-1,0) ..
        (fld-el-1-1) ;

      \draw [thick,blue,->]
        (res) ++(0.3,-0.9)
        .. controls +(0,1) and +(1,0) ..
        (fld-el-2-1) ;
      \draw [thick,blue,->]
        (res) ++(0.3,-1.3)
        .. controls +(1,0) and +(1,0) ..
        (fld-el-2-2) ;

    \end{tikzpicture}
    \caption{
      Row-wise segmentation of a microblocked matrix-matrix product. Element 
      boundaries are shown in black, segment boundaries in red. Also shown:
      Fetch redundancy caused by segmentation. The second segment fetches
      field data from both the first \emph{and} the second element because
      it overlaps rows from both.
    }
    \label{fig:smem-matrix-segment}
  \end{center}
\end{myfigure}

This computation layout makes the shared memory access conflict-free. Unfortunately,
it also introduces a different, smaller drawback: there now is fetch redundancy. 
A segment needs to fetch field data for each element ``touched'' by its rows. 
This may lead it to fetch the same field values as the segment above and below it.  
Figure \ref{fig:smem-matrix-segment} gives an indication of this fetch redundancy,
too.  Fortunately, these duplicated accesses tend to happen in adjacent thread blocks
and therefore possibly at the same time. We speculate that the L2 texture cache in
the device can help reduce the resulting increased bandwidth demand.

Next, observe that the matrix segments typically use less memory than the whole
matrix. We can therefore reexamine the assertion that loading both matrix and
fields into shared memory is not viable. Unfortunately, while the space to do
so is now available, the field access bank conflicts from Section
\ref{ssec:flux-lift} spoil the idea.

One final observation is that for the typical choice of the reference element
\cite{hesthaven_nodal_2007} the three differentiation matrices $D^{\partial
\mu}$ are all similar to each other by a permutation matrix. Using this fact
could allow for significant storage savings, but in our experiments, the added
logic was too costly to make this trick worthwhile.

Algorithm \ref{alg:diff-smem-matrix} presents an overview of the techniques in
this section. Instead of maintaining three separate local differentiation
matrices, it works with one matrix in which the $D^{\partial\mu}$ are
horizontally concatenated and then segmented. Shared memory limitations allow
this algorithm to work at order six and below.

\begin{algorithm}
\caption{Local Differentiation with a segemented matrix in shared memory.}
\label{alg:diff-smem-matrix}
\begin{algorithmic}
  \REQUIRE A grid of $\lceil N_{pM} / N_R \rceil \times \intceil{n_M / (w_p w_i w_s)}{}$ 
    blocks of size $N_R \times w_p \times 1$.
  \REQUIRE Inputs: 
      $\tex u$, the field to be differentiated;
      $\tex r$, the local-to-global differentiation coefficients.
  \ENSURE Output: $\gmem d_\nu$, the local $x,y,z$-derivatives of $\tex u$.
  \STATE Allocate the differentiation matrix chunk
    $\smem{D}\in\mathbb{R}^{N_R\times (N_p d)}$ in shared memory.
  \STATE Load rows $\halfopen{b_x N_R}{b_x (N_R+1)}$ ($\bmod N_p$) 
    of $[D^{\partial 1},\dots, D^{\partial d}]$ 
    into $\smem{D}$.
  \BARRIER
  \FORALL{$\reg s\in[0,w_s)$}
    \STATE $\reg m \leftarrow ((b_yw_s+\reg s)w_p+t_y)w_i$
      \COMMENT{this thread's microblock number}
    \STATE $\reg{d}_\mu^i \leftarrow 0$ 
      for $\mu\in\{1,\dots, d\}$ and $i\in\halfopen{0}{w_i}$
    \FORALL{unrolled $n\in\halfopen 0{N_p}$}
      \STATE $\reg u_{\halfopen{0}{w_i}} \leftarrow 
      \tex u[(\reg m+\halfopen{0}{w_i})N_{pM}+n]$
    \STATE $\reg d_\mu^i \leftarrow \reg d_\mu^i + \smem D[t_x,\mu N_p+ n]\reg u_i$
        for $\mu\in\{1,\dots, d\}$ and $i\in\halfopen{0}{w_i}$
    \ENDFOR
    \STATE $(\gmem{d})_{\halfopen 0 d}^{\reg m N_{pM}+\halfopen{0}{w_i}N_{pM}+t_x} 
      \leftarrow\sum_{\mu} (\tex{r})_{\halfopen{0}{d}d+\mu}^{(m+\halfopen{0}{w_i})K_M} \reg{d}_\mu^i$
  \ENDFOR
\end{algorithmic}
\end{algorithm}

\begin{myfigure}
  \centering
  \tikzset{
    arrowlabeln/.style={
      anchor=south,
      font=\footnotesize,
      text width=1.2cm,
      text centered
    },
    complayoutlabel/.style={
      font=\footnotesize,
      anchor=south,
      rotate=270,
      inner sep=0.02cm
    },
    arrowlabele/.style={
      anchor=west,
      font=\footnotesize
    },
    description/.style={
      anchor=north,font=\footnotesize,
      text width=2cm,
      text centered
    }
  }
  \pgfmathsetmacro{\banks}{16}
  \pgfmathsetmacro{\rowfloats}{20}
  \pgfmathsetmacro{\padfloats}{21}
  \colorlet{evenrow}{black!20}
  \colorlet{oddrow}{black!40}
  \newcommand{\bankcolorize}[2]
  {
    \foreach \row in {#1}
    { 
      \foreach \col in {0,...,19}
      {
        \pgfmathsetmacro{\num}{\row*21+\col}
        \draw [fill=#2] 
          (matbank)
          ++({0.1*mod(\num,\banks)},{0.1*floor(\num/\banks)}) 
          rectangle +(0.1,0.1) ;
      }
    }
  }
  \newcommand\bankedmatrix{
    \draw [step=0.1] (matbank) grid +(1.6,2.7) ;

    \bankcolorize{0,2,...,19}{evenrow}
    \bankcolorize{1,3,...,19}{oddrow}

    \draw (matbank) +(0.8,-0.2) node [description] { Banked Matrix Storage} ;

    \draw [->] (matbank) -- ++(1.8,0) node [arrowlabele] {bank};
    \draw [->] (matbank) -- ++(0,2.9) node [arrowlabeln] {memory};
  }
  \newcommand{\mublock}{
    \draw (mublock) rectangle +(0.75,4.8);
    \draw [fill=gray!20] (mublock) ++(0,0) rectangle ++(0.25,4.0);
    \draw (mublock) ++(0.25,0) -- ++(0,4.8);
    \draw (mublock) ++(0.5,0) -- ++(0,4.8);

    \draw (mublock) ++(0,2) -- ++(0.5,0);
    \draw (mublock) ++(0,4) -- ++(0.5,0);
    \draw (mublock) +(0.5,2.4) node [complayoutlabel] {Microblock} ;
    \draw (mublock) +(0.25,1) node [complayoutlabel] {Element 0} ;
    \draw (mublock) +(0.25,3) node [complayoutlabel] {Element 1} ;

    \draw [->] (mublock) -- ++(0,5) node [arrowlabeln] {thread number};

    \draw [densely dotted,thick, red] 
      (mublock) ++ (0,1.6) -- ($(mublock) + (0.75,1.6) $);
    \draw [densely dotted, thick, red] 
      (mublock) ++(0,3.2) -- ++(0.75,0) ;
    
    \foreach \y in {0.05,0.15,...,3.95}
    { \draw [fill=black] (mublock) ++(0.125,\y) circle (0.015) ; }

    \draw (mublock) +(0.375,-0.2) node [description] { Computation Layout } ;
  }
  \subfigure[
    Assignment of matrix rows to memory banks. 
    Alternating matrix rows are shown in two different shades of gray.
    They preserve their color as they move into individual 4-byte
    cells in the banked shared storage.
    Padding inserted to prevent conflicts is shown in white.
  ]{
    \label{fig:smem-matrix-assignment}
    \begin{tikzpicture}
      \draw [xstep=2,ystep=0.1] (0,0) coordinate (matfull) grid +(2,2) ;

      \coordinate (matbank) at (3,0) ;
      \bankedmatrix

      \newcommand{\matcolorize}[2]{
        \foreach \row in {#1}
        {
          \draw [fill=#2] 
            (matfull)
            ++(0,\row*0.1) 
            rectangle +(2,0.1) ;
        }
      }
      \matcolorize{0,2,...,19}{evenrow}
      \matcolorize{1,3,...,19}{oddrow}

      \foreach \row in {0,1,2,3,19}
      { 
        \draw [blue,->] 
          (matfull) ++(0.05,0.05+0.1*\row)
          .. controls +(2,0.5) and +(-1,1) ..
          ($ (matbank) + (0.05,0.05) + 0.1*({mod(\row*\padfloats,\banks)},{floor(\row*\padfloats/\banks)}) $)
          ; 
      }

      \draw [->] (matfull) -- ++(0,2.2) node [arrowlabeln] {row number};
      \draw (matfull) +(1,-0.2) node [description] { Full Matrix } ;
    \end{tikzpicture}
  }
  \hspace{1em}
  \newcommand{\dofsetup}{
    \pgfmathsetmacro{\dof}{\num*\padfloats}
    \pgfmathsetmacro{\bank}{mod(\dof,\banks)}
    \pgfmathsetmacro{\row}{floor(\dof/\banks)}
  }
  \subfigure[
    Conflict-free access pattern in the first half-warp of the computation layout.
    The green highlighting illustrates that each of the 16 accesses lands in 
    a unique bank.
  ]{
    \label{fig:smem-matrix-conflict-free}
    \begin{tikzpicture}
      \coordinate (matbank) at (0,0) ;
      \bankedmatrix
      \coordinate (mublock) at (3,0) ;
      \mublock

      \foreach \num in {0,...,15}
      {
        \dofsetup
        \draw [fill=darkgreen,opacity=0.5]
          (matbank) ++(0.1*\bank,0) rectangle (0.1*\bank+0.1,0.1*\row+0.1);
      }

      \foreach \num in {0,...,15}
      {
        \dofsetup
        \draw [darkgreen,->]
          ($ (mublock) + (0.05,{0.05+0.1*\num}) $)
          --
          ($ (matbank) 
          + ( {0.05+0.1*\bank}, {0.05+0.1*\row}) $) ;
      }

    \end{tikzpicture}
  }
  \subfigure[
    Conflicting access pattern in the second half-warp of the computation layout.
    The memory banks highlighted in red show 4 banks with two accesses each.
  ]{
    \label{fig:smem-matrix-conflicting}
    \begin{tikzpicture}
      \coordinate (matbank) at (0,0) ;
      \bankedmatrix
      \coordinate (mublock) at (3,0) ;
      \mublock

      \foreach \num in {16,...,19}
      {
        \dofsetup
        \draw [fill=red,opacity=0.5]
          (matbank) ++(0.1*\bank,0) rectangle (0.1*\bank+0.1,0.1*\row+0.1);
      }

      \foreach \num in {16,...,19}
      {
        \dofsetup
        \draw [red,->]
          ($ (mublock) + (0.05,{0.05+0.1*\num}) $)
          --
          ($ (matbank) 
          + ( {0.05+0.1*\bank}, {0.05+0.1*\row}) $) ;
      }

      \foreach \num in {0,...,3}
      {
        \dofsetup
        \draw [red,->]
          ($ (mublock) + (0.05,{2.05+0.1*\num}) $)
          --
          ($ (matbank) 
          + ( {0.05+0.1*\bank}, {0.05+0.1*\row}) $) ;
      }

    \end{tikzpicture}
  }
  \caption{Local matrices and memory banks.}
\end{myfigure}

\section{Experimental Results}

\label{sec:experiments}

In this section, we examine experimental results obtained from a DG solver for
Maxwell's equations in three dimensions for linear, isotropic, and 
time-invariant materials. In terms of the electric field $E$, the magnetic field $H$,
the charge density $\rho$, the current density $j$, the permittivity $\epsilon$,
and the permeability $\mu$, they read
\begin{align}
  \label{eq:maxwells-claws}
  \epsilon \partial_{t}E-\nabla\times H & =-j, 
  & 
  \mu \partial_{t}H+\nabla\times E & =0,\\
  \label{eq:maxwells-elliptic}
  \nabla\cdot (\epsilon E) & =\rho, 
  & 
  \nabla\cdot (\mu H) & =0.
\end{align}
We absorb $E$ and $H$ into a single state vector 
\[
  u\assign (E,H)^T=(E_x, E_y, E_z, H_x, H_y, H_z)^T.
\]
If we define
\[
  F(u)\assign\begin{bmatrix}
    0 & -E_z & E_y & 0 & H_z & -H_y \\
    E_z & 0 & -E_x & -H_z & 0 & H_x \\
    -E_y & E_x & 0 & H_y & -H_x & 0
  \end{bmatrix}^T,
\]
\eqref{eq:maxwells-claws} is equivalently expressed in conservation form as
\[
  \begin{bmatrix} \epsilon & 0 \\ 0 & \mu \end{bmatrix}
  u_t+\nabla\cdot F(u) = 0.
\]
If the two equations \eqref{eq:maxwells-elliptic} are satisfied in the initial
condition, the equations \eqref{eq:maxwells-claws} ensure that this continues
to be the case. Remarkably, the same is also true of the DG discretization of
the operator \cite{hesthaven_nodal_2002}. We may therefore assume a compliant
initial condition and omit \eqref{eq:maxwells-elliptic} from our further
discussion.

We label the numerical solution $u_N\assign (E_N, H_N)^T$ and choose the numerical 
flux $F^*$ to be the upwind flux from \cite{hesthaven_nodal_2002}:
\[
  \hat n \cdot( F_N-F_N^{\ast})\assign
  \frac 12
  \begin{bmatrix}
    \avg Z^{-1} \hat n \times (Z^+ \jump{H_N} - \hat n \times \jump{E_N}) \\
    \avg Y^{-1} \hat n \times (-Y^+ \jump{E_N} - \hat n \times \jump{H_N})
  \end{bmatrix}.
\]
We have employed the conventional notations for the cross-face average 
$\avg u\assign (u_N^-+u_N^+)/2$ and jump $\jump u \assign u_N^+-u_N^-$. For concise
notation, we use the intrinsic impedance $Z\assign \sqrt{\mu/\epsilon}$ and admittance 
$Y\assign 1/Z$. Applying the principles of Section \ref{sec:dg-overview}, we
arrive at a discontinuous Galerkin scheme.

For our experiments, a solver using this scheme runs on an off-the-shelf Nvidia
GTX~280 GPU with 1~GiB of memory using the Nvidia CUDA driver version 177.67. The
GPU code was compiled using the Nvidia CUDA compiler version 2.0.  At the time
of this writing, GPUs of the same type as the one used in this test are sold
for less than US\$400.  

We use a rectangular, perfectly conducting vacuum cavity (see \cite[Section
8.4]{jackson_classical_1998}) excited by one of its eigenmodes to test the
approximate solutions for accuracy. The solver works in single precision.
$L^2$ errors observed for a sequence of grids at orders from one through nine
are shown in Table \ref{tab:convergence}. To better display the actual
convergence of the method, the meshes examined were chosen to be rather coarse.
Between the onset of asymptotic behavior and the saturation at the limits of
single precision, the error exhibits the expected asymptotic behavior of
$h^{N+1}$ \cite{hesthaven_nodal_2002}. We observe that the solver recovers a
significant part of the accuracy provided by IEEE 754 single precision floating
point. It exhibited the same stability properties and CFL time step
restrictions as a corresponding single- and double-precision CPU
implementation.  We have thus established that the discussed algorithm works
and provides solution accuracy on a par to what would be expected of a
single-precision CPU solver. 

\begin{table}[h]
  \centering
  \begin{tabular}{c|cccc|c}
  \hline
  $K$ & 475 & 728 & 1187 & 1844 \\
$N$ & $h=0.3$ & $h=0.255$ & $h=0.21675$ & $h=0.184237$ & EOC \\
\hline
1 & $1.57\cdot10^{0}$ & $1.19\cdot10^{0}$ & $1.03\cdot10^{0}$ & $6.46\cdot10^{-1}$ & 1.72 \\
2 & $4.15\cdot10^{-1}$ & $2.84\cdot10^{-1}$ & $1.82\cdot10^{-1}$ & $1.19\cdot10^{-1}$ & 2.58 \\
3 & $1.61\cdot10^{-1}$ & $9.44\cdot10^{-2}$ & $5.56\cdot10^{-2}$ & $2.80\cdot10^{-2}$ & 3.55 \\
4 & $4.75\cdot10^{-2}$ & $2.52\cdot10^{-2}$ & $1.13\cdot10^{-2}$ & $5.03\cdot10^{-3}$ & 4.64 \\
5 & $1.54\cdot10^{-2}$ & $6.37\cdot10^{-3}$ & $2.55\cdot10^{-3}$ & $9.03\cdot10^{-4}$ & 5.79 \\
6 & $3.84\cdot10^{-3}$ & $1.42\cdot10^{-3}$ & $4.42\cdot10^{-4}$ & $1.32\cdot10^{-4}$ & 6.94 \\
7 & $9.89\cdot10^{-4}$ & $2.77\cdot10^{-4}$ & $7.36\cdot10^{-5}$ & $1.77\cdot10^{-5}$ & 8.24 \\
8 & $1.91\cdot10^{-4}$ & $4.76\cdot10^{-5}$ & $1.05\cdot10^{-5}$ & $2.55\cdot10^{-6}$ & 8.90 \\
9 & $4.25\cdot10^{-5}$ & $8.71\cdot10^{-6}$ & $2.10\cdot10^{-6}$ & $1.30\cdot10^{-6}$ & 7.31 \\

  \hline
  \end{tabular}

  \caption{
    $L^2$ errors and empirical orders of convergence obtained by a solver for
    Maxwell's equations on an Nvidia GTX~280 running in single precision, at a
    variety of orders and for a number of rather coarse meshes.
  }
  \label{tab:convergence}
\end{table}

The reason for bringing DG onto a GPU was however not to show that it works
there, but to show that it can be made to work extremely fast.  Figure
\ref{fig:order-gflops} portrays the speed of our solver in comparison with a
CPU implementation running on a single core of a 3 GHz Intel Core2 Duo E8400
CPU, also running in single precision. The calculations used ATLAS 3.8.2
\cite{whaley_automated_2001} for element-local operations if such use proved
advantageous.  The results are scaled as floating point operations per second,
obtained by counting the number of floating point additions and multiplications
in the algorithm and dividing by the time in seconds.  GPU times were measured
using the \texttt{cuEventElapsedTime()} call. Overall, the GPU outperforms the
CPU by factors ranging from 24 to 57. At the practically relevant orders of
three and four, the speedup factors are 48 and 57, respectively. It is worth
noting that these two orders are not only the ones that see most practical use,
they also exhibit some of the largest speedup factors on the GPU.

\begin{myfigure}
  \centering
  \subfigure[
    Discontinuous Galerkin performance in GFlops/s 
    on a GPU and a CPU. Computations were performed in
    single precision.
  ]{
    \label{fig:order-gflops}
    \includegraphics[width=0.45\textwidth]{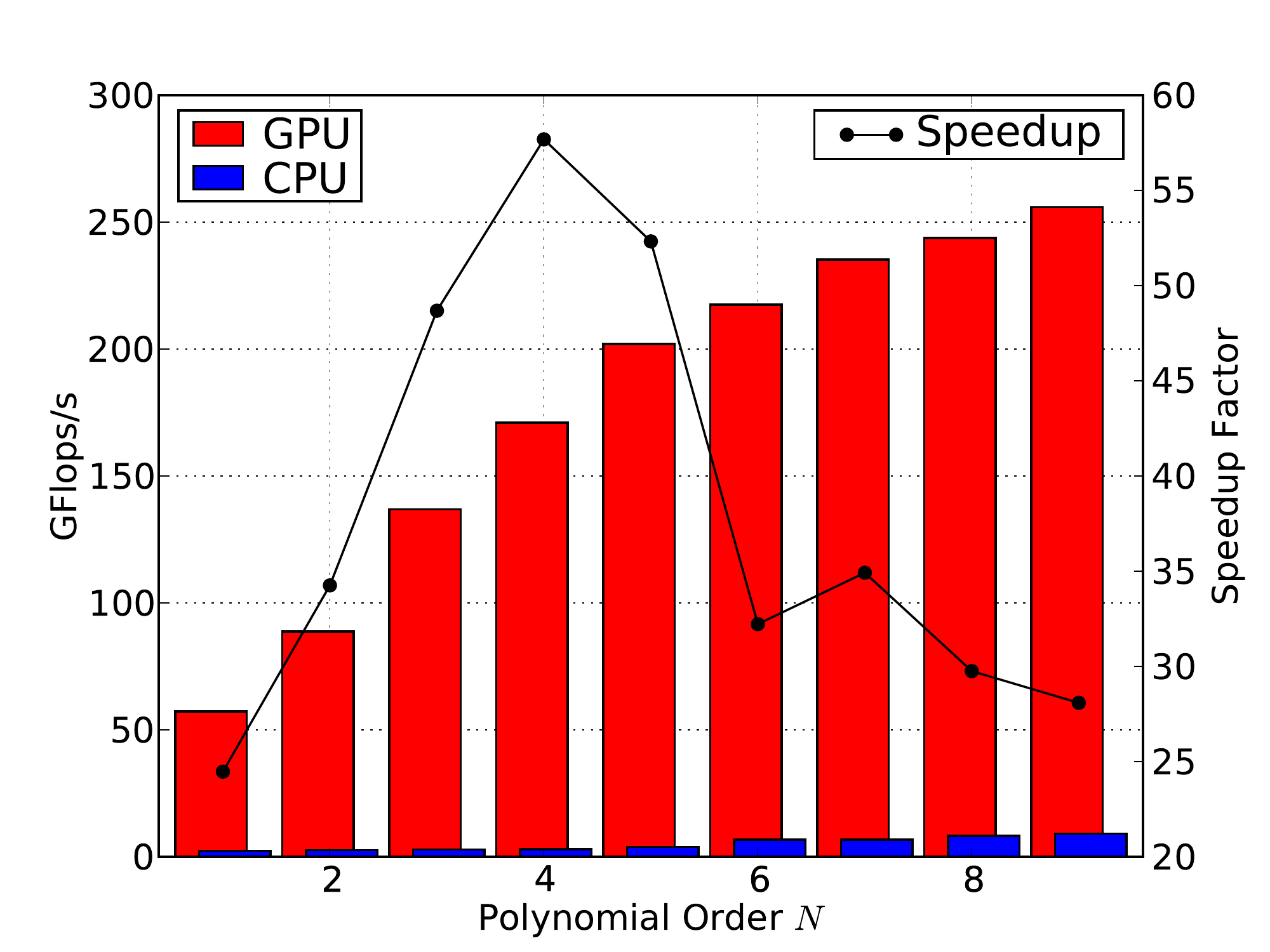}
  }
  \hfill
  \subfigure[
    Number of degrees of freedom to which our methods can apply the
    Maxwelll operator in one second. Assuming linear scaling, this graph
    can be used to determine run times for larger and smaller problems.
    DOFs from each of the six Maxwell fields are counted separately.
  ]{
    \label{fig:dofs-per-sec}
    \includegraphics[width=0.45\textwidth]{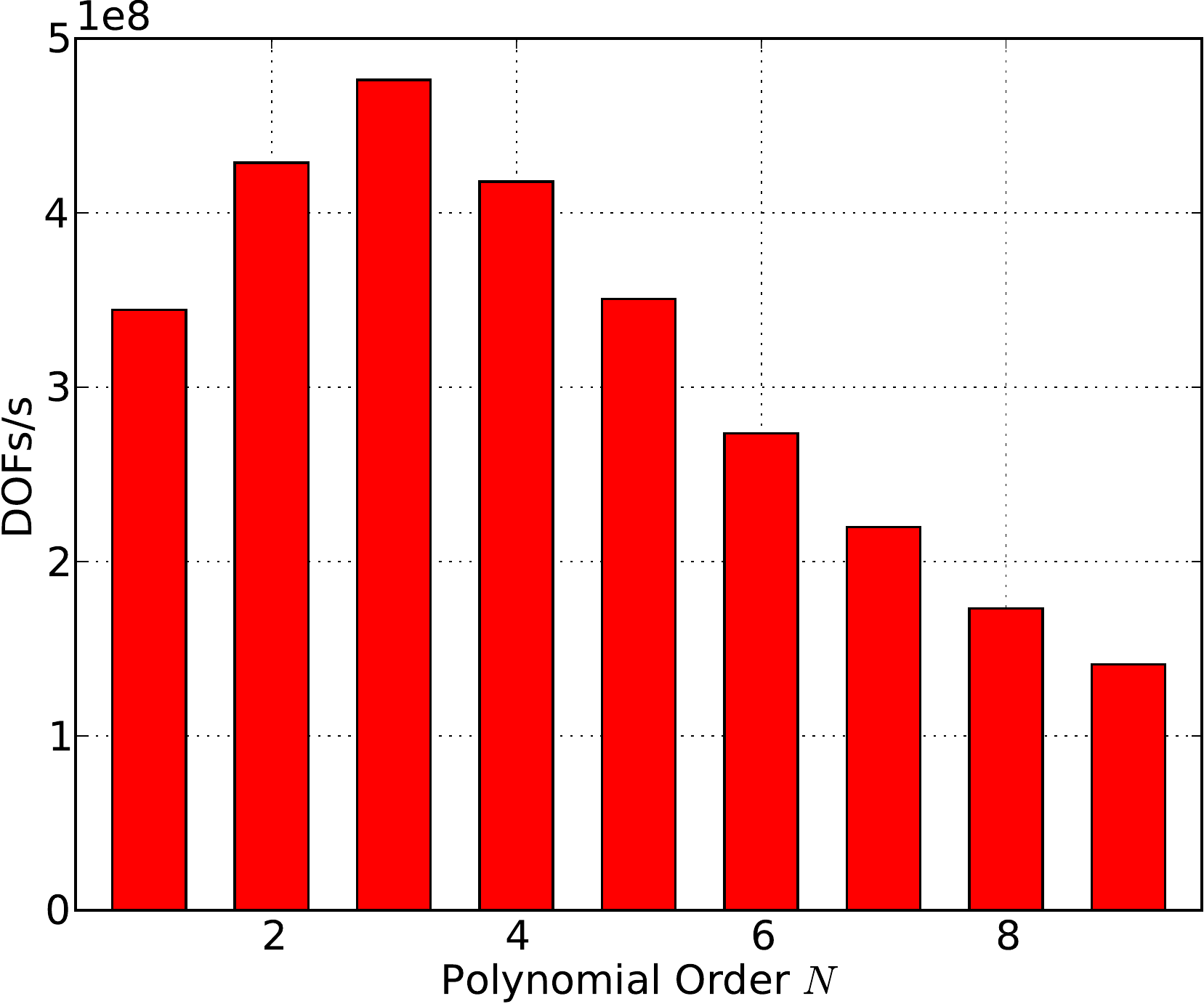}
  }
  \caption{Performance characteristics of DG on Nvidia graphics hardware.}
\end{myfigure}

Orders three and four are particularly favorable not only for their appreciable
speedups and their moderate time step requirements
\cite{warburton_taming_2008}. They also achieve the peak nodal value
throughputs on the GPU as shown in Figure \ref{fig:dofs-per-sec}. Naturally,
high-order approximations of solutions to partial differential equations
contain much more information per DOF than do solutions obtained via low order
methods. This is most apparent in the number of DOFs required to accurately
represent one wavelength \cite{hesthaven_spectral_2007}. Interestingly, we
observe that despite lower computational load, the DG methods of orders one and
two achieve lower overall throughput than the next higher ones, a likely result
of a mismatch with the hardware's granularities.  This crossover between
granularity effects and the increase in floating point work with growing $N$
makes DG methods of orders three and four the fastest DG methods on a GPU even
on a per-DOF basis.

\begin{myfigure}
  \centering
  \subfigure[
    Compute bandwidth in GFlops/s achieved by each part of the DG operator, at
    various polynomial orders.  The published theoretical peak floating point
    performance for the hardware on which these tests were run is 933 GFlops/s
    \cite{wiki:nvidia_comparison}.
  ]{
    \label{fig:gflops-by-component}
    \includegraphics[width=0.45\textwidth]{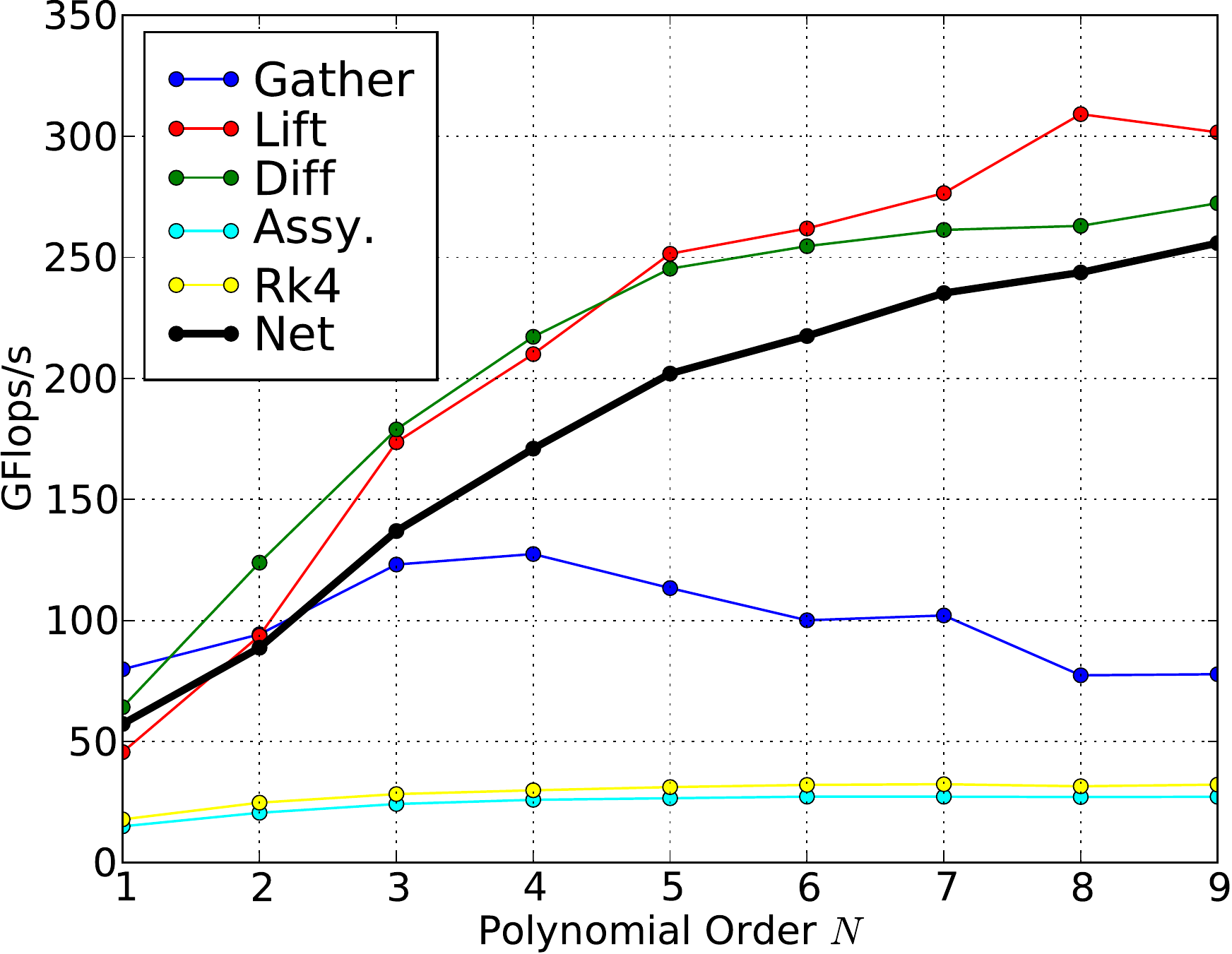}
    }
  \hfill
  \subfigure[
    Percentage of time spent in various parts of the DG operator 
    vs. polynomial order.
  ]{
    \label{fig:time-percentages}
    \includegraphics[width=0.45\textwidth]{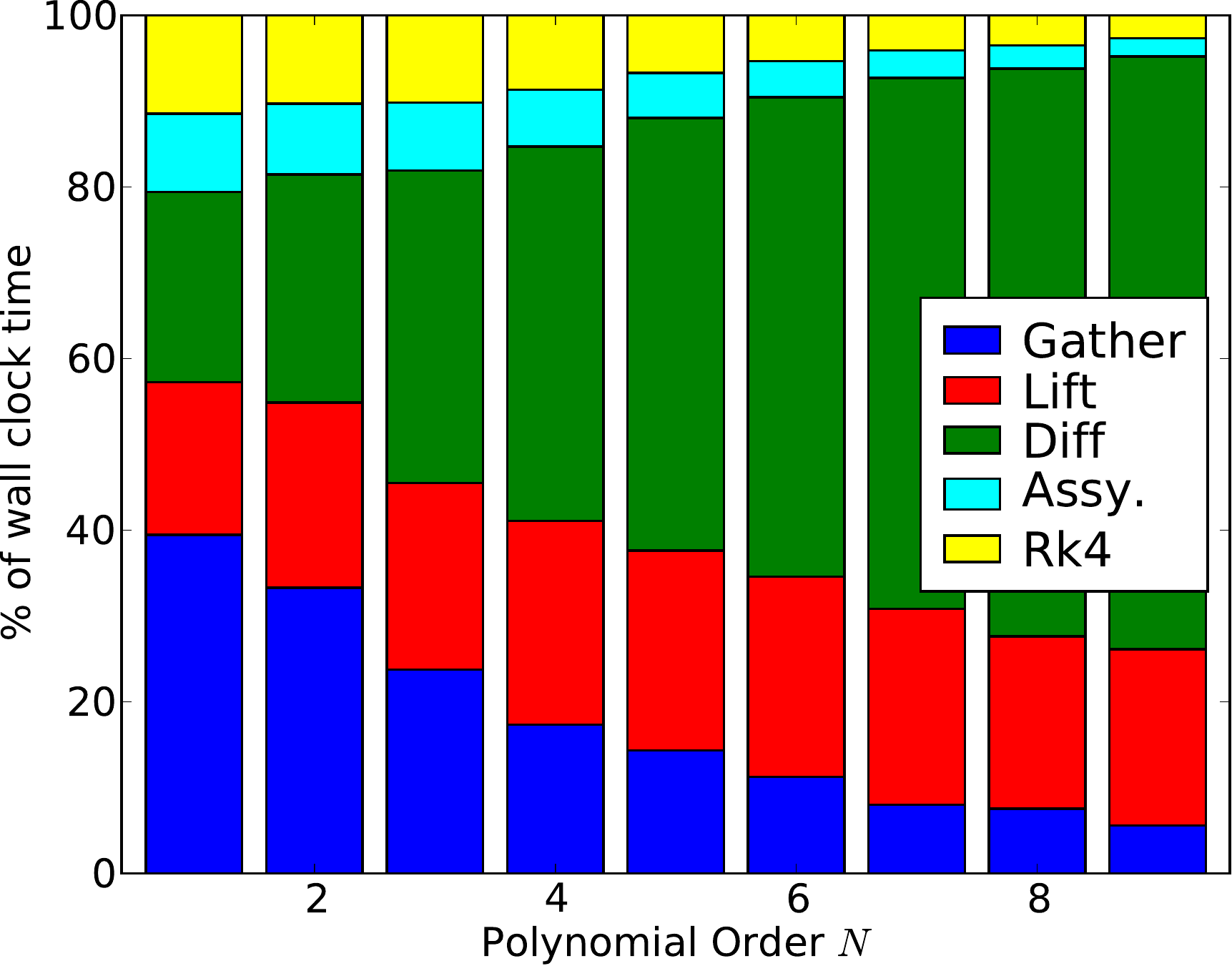}
  }
  \caption{Performance characteristics of DG on Nvidia graphics hardware, continued.}
\end{myfigure}

Recall now that we have split the DG operator into several parts, each of which
performs distinct kinds of processing and, as we have seen, tends to require a
different strategy to map onto a GPU. It is therefore interesting to see what
performance level is attained by each part of the operator. Figure
\ref{fig:gflops-by-component} gives an indication of this performance, based
again on the number of floating point operations per second.  It is reassuring
that, despite different implementation strategies, the flop rates for
element-local differentiation and lifting evolve almost identically as the
order $N$ is increased. These two parts of the operator are also characterized
by the highest arithmetic intensity and the most regular access pattern.  As an
unsurprising consequence, they are able to realize the greatest performance
gain as the order of the operator and therefore the access granularity grows.
The flux gather, on the other hand, realizes its greatest performance at orders
three and four.  We suspect that the decline in performance with increasing $N$
can be attributed to the growth of the indirect indexing information in the
form of face index lists $\tex J$ from Algorithm \ref{alg:flux-extraction}.
These lists are referenced constantly throughout the whole algorithm and are
therefore likely to reside in the texture cache, of which there are only a few
KiB per multiprocessor. As these lists grow, their cache eviction
likelihood also grows, resulting in an increased access latency.  In addition
to the above-mentioned main parts of the operator, the figure also shows
performance data for the assembly of the operator and the fourth-order
low-storage Runge Kutta timestepper \cite{carpenter_fourth-order_1994}. Both of
these operations perform linear combinations of vectors, making them much less
arithmetically intense than the element-local operations.  Fortunately, as the
order $N$ increases, the processing time spent in element-local operations
dominates and helps decrease the influence of the latter three operations on
overall performance. Figure \ref{fig:time-percentages} reinforces this point.

\begin{myfigure}
  \subfigure[
    Memory bandwidths in GB/s achieved by each part of the DG operator.
    The peak memory bandwidth published by the manufacturer is 141.7 GB/s.
    Values exceeding peak bandwidth are believed to be due to the presence
    of a texture cache.
  ]{
    \label{fig:mem-bandwidth}
    \includegraphics[width=0.45\textwidth]{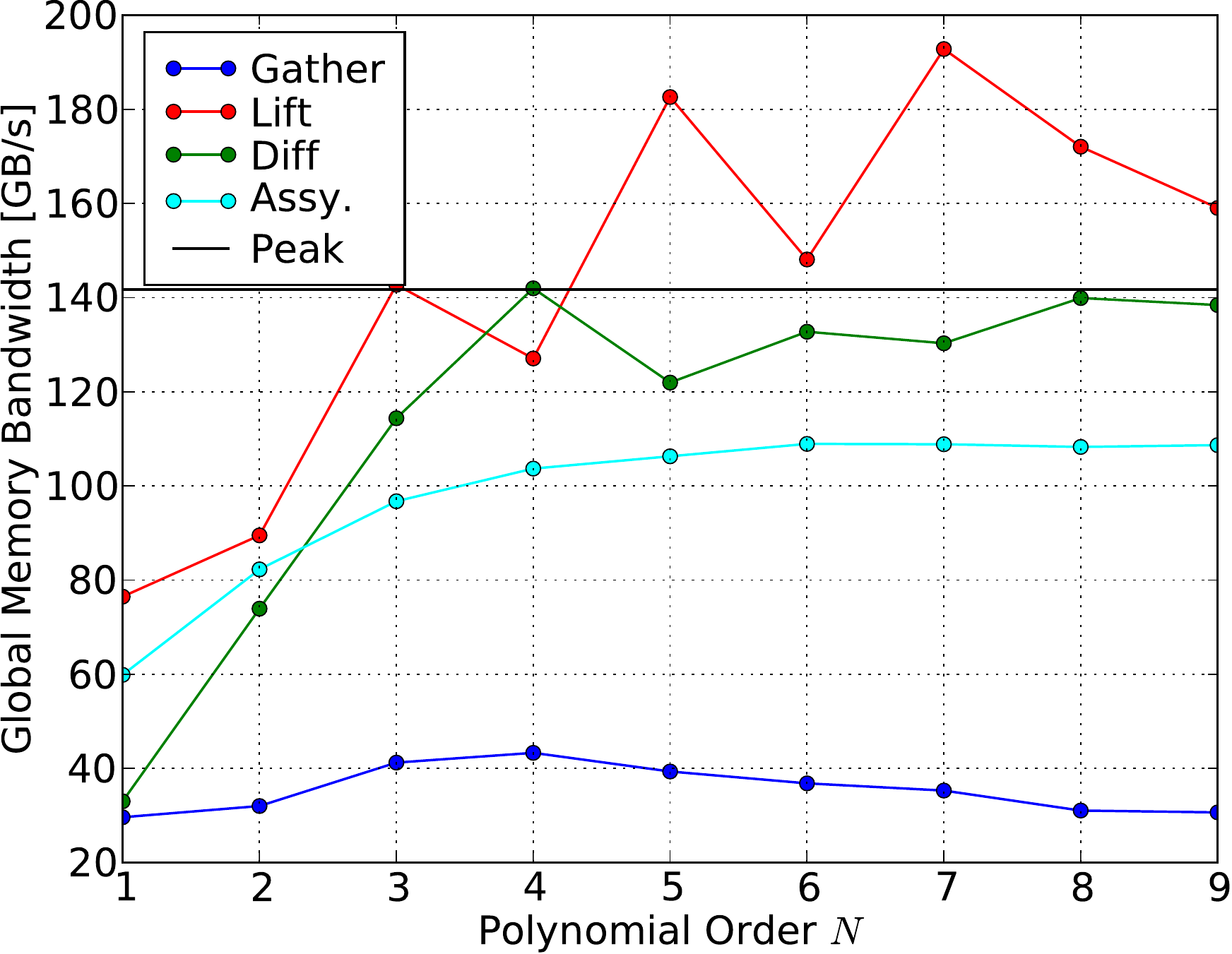}
  }
  \hfill
  \subfigure[
    Sample work distribution parameter study for local differentiation
    on fourth-order elements with microblocking enabled.
  ]{
    \includegraphics[width=0.45\textwidth]{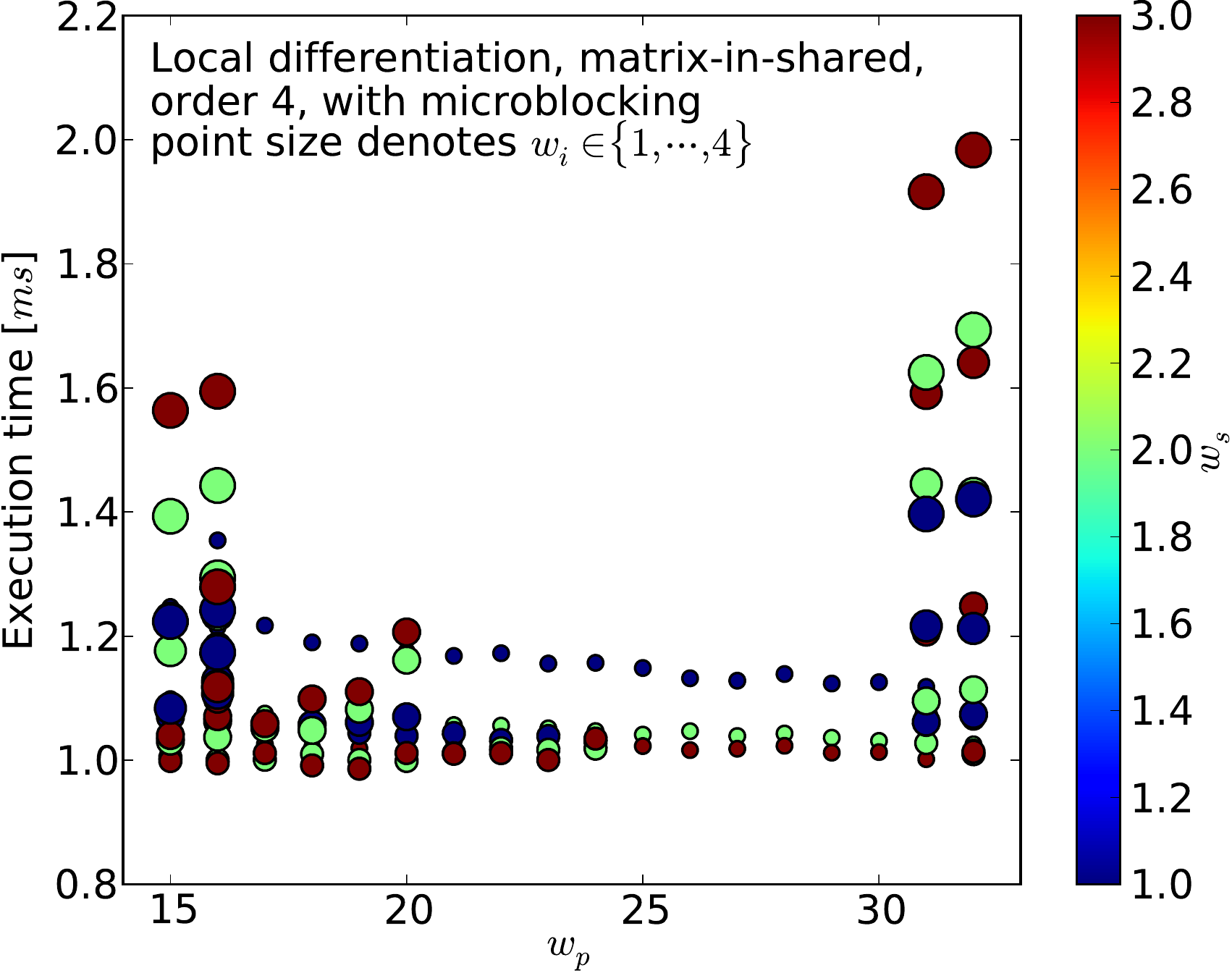}
    \label{fig:diff4-mat-smem-work-distribution}
  }
  \caption{Performance characteristics of DG on Nvidia graphics hardware, continued.}
\end{myfigure}

It is interesting to correlate the achieved floating point bandwidth of each
component from Figure \ref{fig:gflops-by-component} with the bandwidth reached
for transfers between the processing core and global memory, shown in Figure
\ref{fig:mem-bandwidth}. We obtained these numbers by counting the number of
bytes fetched from global memory either directly or through a texture unit. The
theoretical peak memory bandwidth published by Nvidia is 141.7 GB/s, shown as a
black horizontal line. Perhaps the most striking feature at first is the fact
that the memory bandwidth measured for flux lifting transcends this theoretical
peak at orders five and above. We attribute this phenomenon to the presence of
various levels of texture cache. It should perhaps be sobering that the other
parts of the DG operator do not manage the same feat. In any case, flux lifting
uses the fields-in-shared strategy, and therefore fetches and re-fetches the
rather small matrix $L$, making large amounts of data reuse a plausible
proposition. Aside from this surprising behavior of flux lifting, it is both
interesting and encouraging to see how close to peak the memory bandwidth for
element-local differentiation gets. As a converse to the above, this makes it
likely that the operation does not get much use out of the texture cache. It
does imply, however, that rather impressive work was done by Nvidia's hardware
designers: The theoretical peak global memory bandwidth can very nearly be
attained in real-world computations. Next, taking into account what was said in
Section \ref{ssec:flux-lift} about the flux-gather part of the operator, the
rather low memory throughput achieved is not too surprising--the access pattern
is (and, for a general grid, has to be) rather scattered, decreasing the
achievable bandwidth. Lastly, operator assembly, which computes linear
combination of vectors, consists mainly of global memory fetches and stores.
There is no reason why it should be unable to pin the memory subsystem to its
peak throughput.  Unfortunately, we found ourselves unable to achieve this
through loop unrolling or other techniques.

\begin{table}
\begin{center}
\begin{tabular}{c|c|cccc|cc|cccc}
\hline
 & & \multicolumn{4}{|c}{differentiation}
 & \multicolumn{2}{|c}{flux gather} 
 & \multicolumn{4}{|c}{flux lifting} 
 \\
$N$ & $K_M$ 
& Shared & $w_p$ & $w_i$ & $w_s$ 
& $M_B$ & $w_p$ 
& Shared & $w_p$ & $w_i$ & $w_s$
 \\
\hline
1 & 4 & Matrix & 6 & 2 & 3 & 4 & 28 & Field & 4 & 2 & 1 \\
2 & 8 & Matrix & 19 & 1 & 3 & 2 & 26 & Field & 3 & 3 & 1 \\
3 & 4 & Matrix & 14 & 2 & 3 & 3 & 19 & Field & 2 & 3 & 1 \\
4 & 4 & Matrix & 19 & 2 & 3 & 2 & 18 & Field & 2 & 4 & 1 \\
5 & 2 & Field & 1 & 4 & 1 & 3 & 15 & Field & 2 & 3 & 1 \\
6 & 2 & Field & 1 & 4 & 1 & 1 & 4 & Field & 2 & 4 & 1 \\
7 & 2 & Field & 2 & 4 & 1 & 2 & 8 & Field & 2 & 3 & 1 \\
8 & 1 & Field & 2 & 4 & 1 & 1 & 1 & Field & 2 & 4 & 1 \\
9 & 1 & Field & 2 & 4 & 1 & 1 & 2 & Field & 2 & 4 & 1 \\

\hline
\end{tabular}
\end{center}
\caption{Empirically optimal method parameters for each part of the DG operator 
  at polynomial orders 1 through 9.}
\label{tab:method-parameters}
\end{table}

For potential implementers, it may be interesting to know which exact
parameters were used to obtain the results in this section.  The parameters of
interest include the generic work distribution tuple $(w_p, w_i, w_s)$ for each
subtask, the microblock size $K_M$, the gather block size $M_B$, and which of
the matrix- or field-in-shared approaches was used at what order. Table
\ref{tab:method-parameters} presents this data. It is peculiar how little
regularity there is in this data set. Despite a sequence of attempts, we failed
to come up with a heuristic that would predict performance accurately. This led
us to develop an empirical optimization procedure that finds the data of Table
\ref{tab:method-parameters} in an automated fashion through a sequence of
synthetic and real-world benchmarks.  A detailed study of this and other
optimization procedures as well as of the toolkit we constructed to enable them
will be the subject of a forthcoming report. For now, we restrict ourselves to
displaying the results of one such procedure. Figure
\ref{fig:diff4-mat-smem-work-distribution} displays the run time obtained for
element-local differentiation employing microblocking and the matrix-in-shared
strategy at order $N=4$. The objective is to find the work distribution
parameter tuple $(w_p, w_i, w_s)$ that leads to an empirically short run time
for this part of the operator. It should be stressed that all runs depicted in
the figure perform the same amount of work. From Table
\ref{tab:method-parameters} we see that in this particular instance, an
optimum was found at $(w_p,w_i,w_s)=(19,2,3)$. Undoubtedly, with better
knowledge of the hardware, many of the odd-looking ups and downs in Figure
\ref{fig:diff4-mat-smem-work-distribution} could be understood. Given the
published documentation however, we are mostly left to take the results at face
value. Luckily, if one were to randomly choose a configuration from the
portrayed set, in all likelihood the resulting operation would at most take
about 20 per cent longer than the optimal one chosen here.  On the other hand,
with some bad luck one may also encounter a configuration that makes the
computation take about twice as long.

\begin{myfigure}
  \centering
  \subfigure[
    Performance in GFlops/s achieved at various polynomial orders,
    for different simplified implementations of DG on CUDA.
  ]{
    \label{fig:gflops-simplified}
    \includegraphics[width=0.45\textwidth]{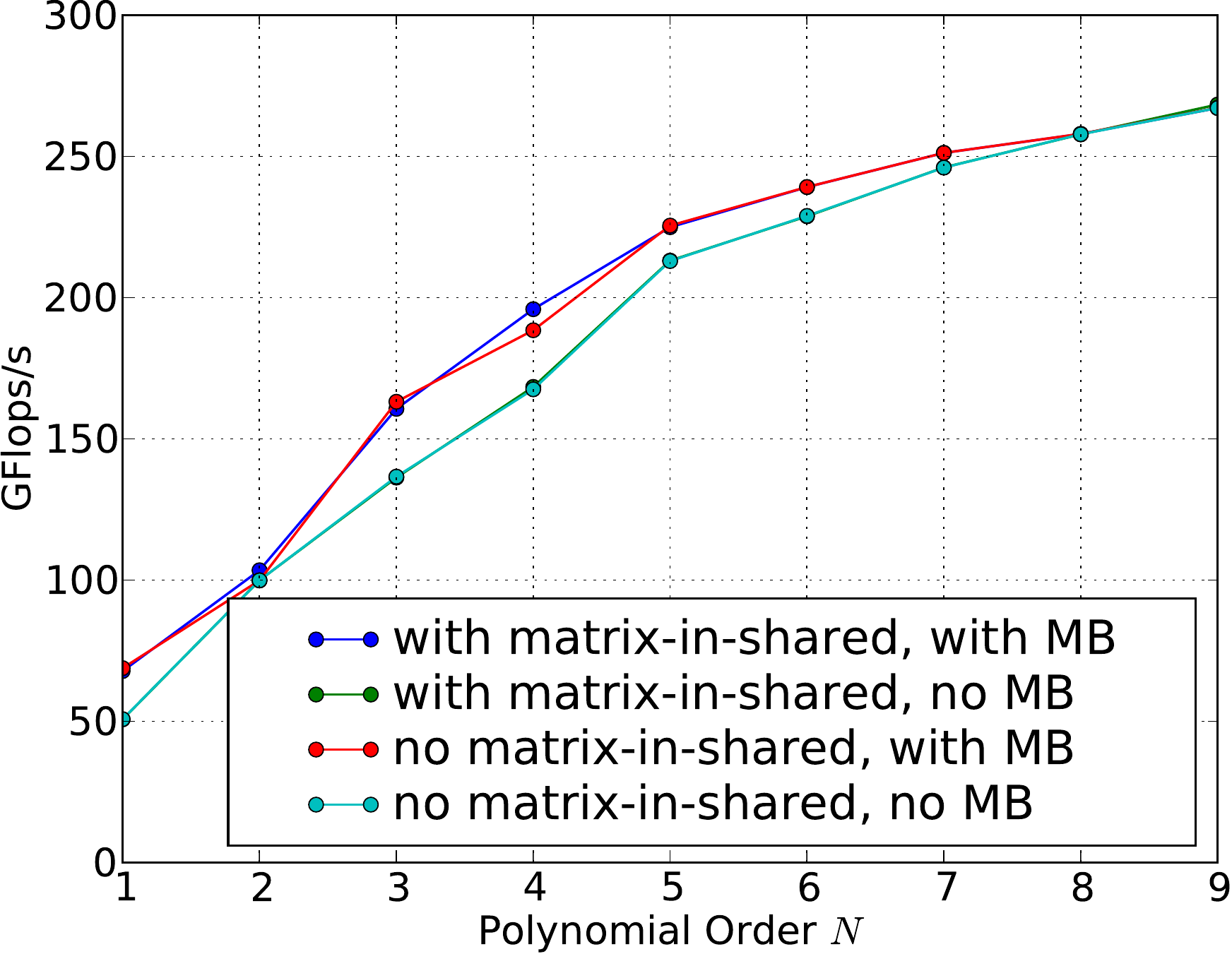}
  }
  \hfill
  \subfigure[
    Mesh-dependent scaling of discontinuous Galerkin on Nvidia GPUs.
  ]{
    \label{fig:gflops-vs-mesh}
    \includegraphics[width=0.45\textwidth]{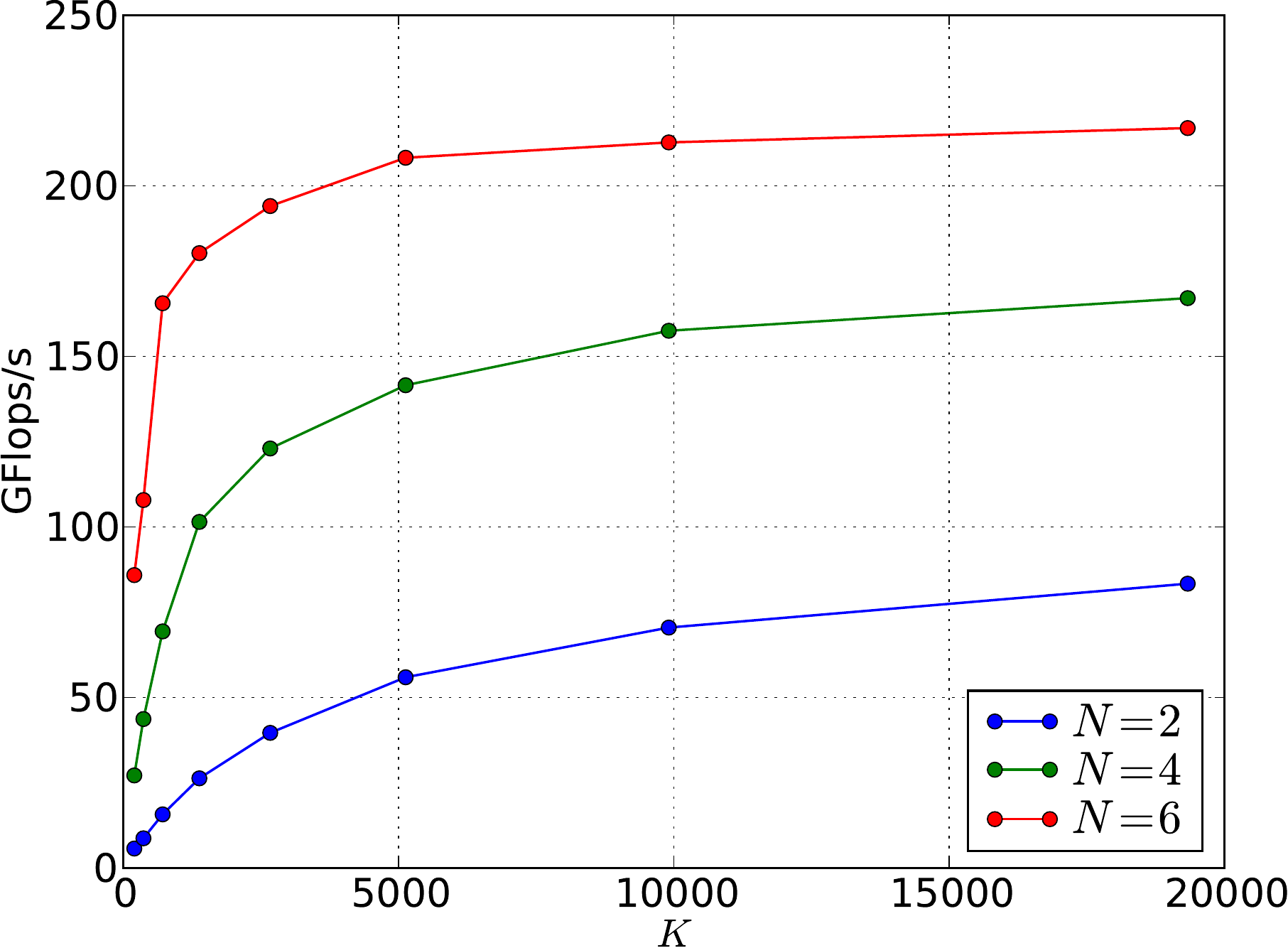}
  }
  \caption{Performance characteristics of DG on Nvidia graphics hardware, continued.}
\end{myfigure}

From Table \ref{tab:method-parameters} we can also gather that the
field-in-shared strategy with a wide variety of work distribution parameters is
found to deliver the best performance at all orders for flux lifting, as well
as for higher-order element-local differentiation. This is plausible behavior
and was already discussed in Section \ref{ssec:local-diff}. It is therefore
reasonable to ask what would be lost if the matrix-in-shared approach were
omitted from a GPU DG implementation entirely. Also, we have seen in a number
of sections that the introduction of microblocks into the method brings about
some mild complications, particularly in the form of shared memory bank
conflicts, so one may be compelled to ask how much is lost by ignoring
microblocks and simply padding each element to the nearest alignment boundary.
The remaining performance after restricting our implementation to not use one
or both of these optimizations can be seen in Figure
\ref{fig:gflops-simplified}.  Examination of this figure leads to the
conclusion that the work of implementing a matrix-in-shared strategy is likely
only worthwhile if one is particularly interested in running GPU-DG at a few
specific low orders. The benefit of employing mircoblocking, on the other hand,
is pervasive and fairly substantial. It stretches to far higher orders than one
might suspect at first, given the growth of the involved operands.  

Note that these conclusions apply only to the algorithms exactly as described
so far.  If even one simple trick is omitted from an implementation, tradeoffs
may shift dramatically. For example, omitting the thread ordering trick from
Section \ref{ssec:flux-lift} makes a matrix-in-shared strategy optimal for
differentiation up to order six.

Finally, we note that the performance results in this section depend on the
size of the problem being worked on. A very small problem may, for example, not
offer enough opportunity to properly occupy all the processing cores that the
hardware provides. Figure \ref{fig:gflops-vs-mesh} reveals that even relatively
small problems achieve decent performance.  In addition, we observe that this
scaling effect is apparently not just governed by the number of elements
present, but also by the order $N$, which influences the number of flops per
DOF in the method. We conclude that as soon as there is a certain amount of
floating point work to be done per timestep, the method will perform fine.

\begin{figure}
  \centering
  \includegraphics[width=0.5\textwidth]{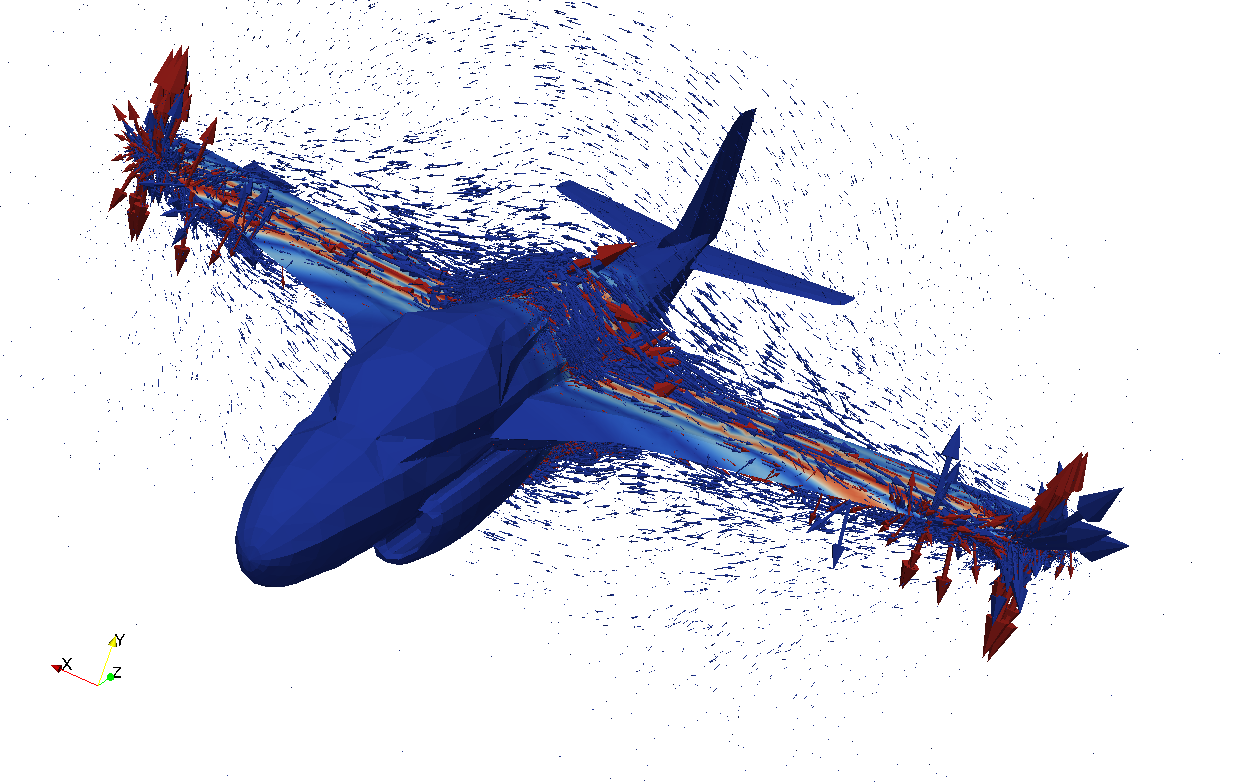}
  \label{fig:sample-calculation}
  \caption{
    A sample scattering problem solved using the methods described in the text.
    The incident plane-wave electric field is shown as pseudocolor values on
    the scatterer, while the scattered electric field is shown as arrows.  The
    computation was performed at order $N=4$ on a mesh of $K=78745$ elements
    using an incident-field formulation \cite{hesthaven_nodal_2002} and
    characteristic absorbing boundary conditions. It achieved and sustained
    more than 160 GFlops/s.
    }

\end{figure}

\section{Conclusions}

\label{sec:conclusions}

In this paper, we have described and evaluated a variety of techniques for
performing discontinuous Galerkin simulations on recent Nvidia graphics
processors.  We have adapted a number of pre-existing DG codes for the GPU,
enabling a thorough comparison of strategies for mapping the method onto the
hardware. A final code was written that combines the insights gained from its
precursors.  This code implements the strategies of Sections
\ref{sec:dg-gpu-design} and \ref{sec:dg-gpu-implementation} and was used to
obtain the results in Section \ref{sec:experiments}.

We have shown that, using our strategies, high-order DG methods can reach
double-digit percentages of published theoretical peak performance values for
the hardware under consideration. DG computations on GPUs see speed-up factors
just short of two orders of magnitude. This speed increase translates directly
into an increase of the size of the problem that can be treated using these
methods.  A single compute device can now do work that previously required a
roomful of computing hardware. Alternatively, a cluster of machines equipped
with these cards can run simulations that were previously outside the reach of
all but the largest supercomputers.  This lets the size and complexity of
simulations that researchers can afford on a given hardware budget jump
significantly.

To make these speed gains accessible, we have provided detailed advice for
potential implementers who need to achieve a trade-off between computing
performance and implementation effort. The data provided in Section
\ref{sec:experiments} will help them make informed implementation decisions by
allowing them to predict the computational speed achieved by their
implementations.

We will be extending this work to make use of double precision floating point
support that has become available on recent Nvidia hardware. In addition, we
would like to broaden the applicability of our methods by exploring their use
for nonlinear conservation laws as well as elliptic problems.

Many-core computing presents a rare opportunity, and we feel that discontinuous
Galerkin methods have a number of unique characteristics that make them
unusually suitable for many-core platforms. In the past, users have chosen
low-order methods because of the perceived expense involved in running
simulations at a high order of accuracy. While this perception was questionable
even in the past, we feel that many-core architectures disproportionately
\emph{favor} high order and significantly drive down its relative cost.
Moreover, unlike most other numerical schemes for solving partial differential
equations, DG methods make the order of accuracy a tunable parameter. These
factors combine to give the user a maximum of control over both performance
and accuracy.

\subsection{Acknowledgments}

The authors gratefully acknowledge the support of AFOSR under grant number
FA9550-05-1-0473.  The opinions expressed are the views of the authors. They do
not necessarily reflect the official position of the funding agencies.

We would like to thank Nvidia Corporation, who, upon completion of
this work, provided us with a generous hardware donation for further
research. 

We would also like to thank Nico G\"odel, Akil Narayan, and Lucas Wilcox who
provided helpful insights in numerous discussions.

Meshes used in this work were generated using TetGen \cite{si_meshing_2005}.
The surface mesh for Figure \ref{fig:sample-calculation} originates in the
FlightGear flight simulator and was processed using Blender and MeshLab, 
a tool developed with the support of the Epoch European Network of Excellence.

\appendix

\section{Index of Notation}
\label{sec:notation}

\begin{tabular}{lp{7.5cm}r}
\hline
Symbol & Meaning & See \\
\hline
$\intceil x n$ & $x$ rounded up to the nearest multiple of $n$. 
  & \ref{sec:pseudocode-notation} \\
$\halfopen a b$ & The set of integers from the half-open interval $[a,b)$. 
  & \ref{sec:pseudocode-notation}  \\
$d$ & The number of dimensions. & \ref{sec:dg-overview} \\
$n$ & The number of unknowns in the conservation law \eqref{eq:claw}. 
  & \ref{sec:dg-gpu-design} \\
$N$ & Polynomial degree of the approximation space. & \ref{sec:dg-overview} \\
$N_p$ & Number of modes/points in local expansion.& \ref{sec:dg-gpu-design} \\
$N_{fp}$ & Number of facial nodes in reference element.& \ref{sec:dg-gpu-design} \\
$N_f$ & Number of faces in the reference element.& \ref{sec:dg-gpu-design} \\
$k$ & Used to refer to element numbers. & \ref{sec:dg-overview} \\
$K$ & Total number of elements. & \ref{sec:dg-overview} \\
$\D_k$ & The $k$th finite element. & \ref{sec:dg-overview} \\
$\mathsf{I}$ & The unit finite element. & \ref{ssec:implement-dg} \\
$\Psi_k$ & The local-to-global map for element $k$. & \ref{ssec:implement-dg} \\
$M^k$, $M^{k,A}$, $L^k$ 
  & Global mass,  face mass and lifting matrices for element $k$.
  & \ref{sec:dg-overview} \\
$S^{k,\partial \nu}$, $D^{k, \partial \nu}$ 
  & $\nu$th global stiffness and differentiation matrices.
  & \ref{sec:dg-overview} \\
$M$, $M^A$, $L$ 
  & Reference mass,  face mass and lifting matrices.
  & \ref{ssec:implement-dg} \\
$S^{\partial \mu}$, $D^{\partial \mu}$ 
  & $\mu$th reference stiffness and differentiation matrices.
  & \ref{ssec:implement-dg} \\
$\nu$ & Used to index global derivatives. 
  & \ref{ssec:implement-dg} \\
$\mu$ & Used to index local derivatives. 
  & \ref{ssec:implement-dg} \\
%
%
$T$ & Thread scheduling (``warp'') granularity. & \ref{sec:gpu-hardware} \\
$K_M$ & Number of elements in one microblock. & \ref{sec:dg-gpu-design} \\
$N_{pM}$ & Number of volume DOFs in a microblock after padding. & \ref{sec:dg-gpu-design} \\
$N_{fM}$ & Number of face DOFs in a microblock after padding. & \ref{ssec:flux-lift}\\
$M_B$ & Number of microblocks in one flux-gather block. 
  & \ref{ssec:flux-gather} \\
$n_M$ & Total number of microblocks. ($=\lceil K/K_M\rceil$) & \ref{sec:dg-gpu-design} \\
$N_R$ & Row count of a matrix segment. & \ref{ssec:local-diff} \\
$w_p$ & The number of work units one block processes in parallel, 
  in different threads. & \ref{sec:dg-gpu-design} \\
$w_i$ & The number of work units one block processes inline, 
  in the same thread. & \ref{sec:dg-gpu-design} \\
$w_s$ & The number of work units one block processes sequentially, 
  in the same thread. & \ref{sec:dg-gpu-design} \\
$t_x,t_y,t_z$ & Thread indices in a thread block.
  & \ref{sec:gpu-hardware} \\
$b_x,b_y$ 
  & Block indices in an execution grid. & \ref{sec:gpu-hardware} \\
\hline
\end{tabular}

\bibliography{dg-cuda}
\bibliographystyle{plainnat}
\end{document}